\documentclass[11pt]{article} 
\usepackage[sectionbib]{natbib}
\usepackage{array,epsfig,fancyheadings,rotating}
\usepackage[]{hyperref}  
\usepackage{sectsty, secdot}
\sectionfont{\fontsize{12}{14pt plus.8pt minus .6pt}\selectfont}
\subsectionfont{\fontsize{12}{14pt plus.8pt minus .6pt}\selectfont}

\usepackage{amsmath}
\usepackage{amssymb}
\usepackage{amsfonts}
\usepackage{multirow}
\usepackage{amsthm}

\newtheorem{theorem}{Theorem}
\newtheorem{lemma}{Lemma}
\newtheorem{corollary}{Corollary}

\theoremstyle{definition}

\newtheorem{remark}{Remark}

\usepackage[utf8]{inputenc}
\usepackage{lipsum}

\usepackage[dvipsnames]{xcolor}
\usepackage{epsfig}
\usepackage{graphicx}
\usepackage{epstopdf}
\usepackage{multirow}
\usepackage{listings}

\graphicspath{{./}}

\usepackage{xr}

\usepackage{float}
\allowdisplaybreaks[1]

\numberwithin{equation}{section}
\numberwithin{theorem}{section}
\numberwithin{figure}{section}
\numberwithin{table}{section}
\numberwithin{remark}{section}
\numberwithin{corollary}{section}

\usepackage{a4wide}

\begin{document}



\markright{ \hbox{\footnotesize\rm 
}\hfill\\[-13pt]
\hbox{\footnotesize\rm
}\hfill }

\markboth{\hfill{\footnotesize\rm W. Li, Q. Wang and J. Yao.} \hfill}
{\hfill {\footnotesize\rm High dimensional distance covariance matrix} \hfill}

\renewcommand{\thefootnote}{}
$\ $\par


\fontsize{12}{14pt plus.8pt minus .6pt}\selectfont \vspace{0.8pc}
\centerline{\large\bf Eigenvalue distribution of a high-dimensional }
\vspace{2pt} 
\centerline{\large\bf distance  covariance matrix with application}
\vspace{.4cm} 
\centerline{Weiming Li, Qinwen Wang and Jianfeng Yao} 
\vspace{.4cm} 
\centerline{\it Shanghai University of Finance and    Economics,}
\centerline{\it Fudan University and The University of Hong Kong}

\vspace{.55cm} \fontsize{9}{11.5pt plus.8pt minus.6pt}\selectfont


\begin{quotation}
\noindent {\it Abstract:} \qquad
  We introduce a new random matrix model called  distance covariance matrix in this paper, whose normalized trace is equivalent to the distance covariance. We first derive a deterministic limit for the eigenvalue distribution of the distance covariance matrix when the dimensions of the vectors and the sample size tend to infinity simultaneously. This limit is valid when the vectors are independent or weakly dependent through a finite-rank perturbation. It is also  universal and independent of the details of the distributions of the vectors. Furthermore, the top eigenvalues of this distance covariance matrix are shown to obey an exact phase transition when the dependence of the vectors is of finite rank. This finding enables the construction of a new detector for such weak dependence where classical methods based on large sample covariance matrices or sample canonical correlations may fail in the considered high-dimensional framework.

\vspace{9pt}
\par
\end{quotation}\par

\fontsize{12}{14pt plus.8pt minus .6pt}\selectfont

%
\newcommand{\1}{{\bf 1}}
\newcommand{\bx}{{\bf x}}
\newcommand{\by}{{\bf y}}
\newcommand\red[1]{{#1}}
\newcommand\bl[1]{{\color{blue}#1}}
\newcommand{\ba}{{\mathbf a}}
\newcommand{\bb}{{\mathbf b}}
\newcommand{\bu}{{\mathbf u}}
\newcommand{\X}{{\mathbf X}}
\newcommand{\Y}{{\mathbf Y}}
\newcommand{\Z}{{\mathbf Z}}
\newcommand{\I}{{\mathbf I}}
\newcommand{\T}{{\mathbf T}}
\newcommand{\A}{{\mathbf A}}
\newcommand{\B}{{\mathbf B}}
\newcommand{\C}{{\mathbf C}}
\newcommand{\D}{{\mathbf D}}
\newcommand{\V}{{\mathbf V}}
\newcommand{\M}{{\mathbf M}}
\newcommand{\bS}{{\mathbf S}}
\newcommand{\W}{{\mathbf W}}
\newcommand{\s}{{\mathbf s}}
\newcommand{\w}{{\mathbf w}}
\newcommand{\e}{{\mathbf e}}
\newcommand{\xx}{{\mathbf x}}
\newcommand{\y}{{\mathbf y}}
\newcommand{\z}{{\mathbf z}}
\newcommand{\br}{{\mathbf r}}
\newcommand{\bP}{{\mathbf P}}
\newcommand{\zt}{{\tilde z}}
\newcommand{\bv}{{\mathbf v}}
\newcommand{\bd}{{\mathbf d}}
\newcommand\tx{\widetilde{\bx}}
\newcommand{\al}{{\alpha}}
\newcommand{\te}{{\theta}}
\newcommand{\tetan}{{\hat\theta_n}}
\newcommand\tetanmk{{\hat{\theta}_{n}^{(k)}}}
\newcommand{\CC}{\mathbb{C}}
\newcommand{\la}{\lambda}
\newcommand{\si}{\sigma}
\newcommand{\Sig}{{\bf \Sigma}}
\newcommand{\Ps}{\boldsymbol \psi}

\newcommand{\um}{\underline{m}}
\newcommand{\uV}{\underline{V}}
\newcommand{\uS}{\underline{S}}
\newcommand{\uT}{\underline{T}}
\newcommand{\de}{\delta}
\newcommand{\E}{{\mathbb E}}
\newcommand{\var}{{\mathbb V}ar}

\newcommand{\Cov}{{\rm Cov}}
\newcommand{\bL}{{\mathbf L}}
\newcommand{\U}{{\mathbf U}}
\newcommand{\bV}{{\mathbf V}}
\newcommand{\bM}{{\mathbf M}}
\newcommand{\bR}{{\mathbf R}}
\newcommand{\bQ}{{\mathbf Q}}
\newcommand{\R}{{\tilde {\mathbf R}}}
\newcommand{\tr}{{\text{\rm tr}}}
\newcommand{\btheta}{{\boldsymbol \theta}}
\newcommand{\be}{{\boldsymbol \varepsilon}}
\newcommand{\bbeta}{{\boldsymbol \beta}}
\newcommand{\bgamma}{{\boldsymbol \gamma}}
\newcommand{\no}{\nonumber\\}
\newcommand{\overbar}[1]{\mkern 1.5mu\overline{\mkern-1.5mu#1\mkern-1.5mu}\mkern 1.5mu}

\newcommand{\mV}{{\mathcal V}}

\newcommand\gai{\color{blue}}

\section{Introduction}\label{sec:intro}

\cite{S07} introduced the concept of  {\em distance covariance}  ${\mathcal V}(\bf x, \bf y)$ of  two random vectors $(\bx, \by)\in\mathbb R^p\times \mathbb R^q$ as a measure of their dependence. 
It is defined through an  appropriately   weighted $L_2$-distance between the joint characteristic function
$\phi_{\bf x,  \bf y}(s,t)$ of $(\bf x,\bf y)$  and the product of their marginal characteristic functions $\phi_{\bf x} (s)  \phi_{\bf y}(t)$,
namely                         
\begin{align}\label{eq:V2}
  {\mathcal V}(\bf x,\bf y)&   
  =\left\{
  \frac{1}{c_p
    c_q}\int\!\!\!\int_{\mathbb{R}^{p}\times\mathbb{R}^{q}}\frac{|\phi_{\bf
      x, \bf y}(s,t)-\phi_{\bf x}(s)\phi_{\bf
      y}(t)|^2}{\|s\|^{1+p}\|t\|^{1+q}}ds dt
  \right\}^{\frac12},  
\end{align}
where the normalization constants are  $c_d = \pi^{(1+d)/2} /\Gamma((1+d)/2)$ $(d=p,q)$.
Clearly, ${{\mathcal V}(\bf x, \bf y)}=0$
if and only if $\bf x$ and $\bf y$ are independent.

{
For a collection of i.i.d.\ observations $(\bx_1, \by_1), \ldots, (\bx_n, \by_n)$ from the population $(\bx,\by)$,
\citet{S07}  proposed   the {\em sample  distance covariance} $\mathcal{V}_n( \bx,\by)$ as 
\begin{equation}
  \label{eq:Vn}
  \mathcal{V}_n( \bx,\by)= \left\{S_{1,n} + S_{2,n}- 2S_{3,n}\right\}^{1/2},
\end{equation}
where
 \begin{align*}
  \nonumber
  S_{1,n}&=\frac{1}{n^2}\sum_{k, \ell=1}^n \|\bx_k-\bx_\ell\| \|\by_k-\by_\ell\|,\quad\\
  S_{2,n}&=\frac{1}{n^2}\sum_{k, \ell=1}^n \|\bx_k-\bx_\ell\| \frac{1}{n^2}\sum_{k, \ell=1}^n \|\by_k-\by_\ell\|,\\
  S_{3,n}& =\frac{1}{n^3}\sum_{k, \ell,m=1}^n \|\bx_k-\bx_\ell\| \|\by_k-\by_m\|.
 \end{align*}
}
One remarkable result \citep[Theorem 2]{S07} is that whenever
$\E[\|\bx\|+\|\by\|]<\infty$, 
$\mV_n(\bx,\by)$ converges almost surely to  $\mathcal{V}( \bx,\by)$
as  $n\to\infty$.
Based on this,   a
powerful statistic
\begin{align}\label{tnsta}
 T_n = n\mV^2_n(\bx,\by)/  S_{2,n}
 \end{align}
  was developed  for  testing the
 independence hypothesis, 
\begin{align}\label{hoo}
H_0: \bx~\text{is~independent~of~}\by,
\end{align}
by establishing:
(i) under  $H_0$,  $\displaystyle  T_n  \overset{\mathcal
  D}{\longrightarrow}  Q$, a countable  mixture of independent
chi-squared distributions,  
and (ii)
  if $\bx$ and $\by$ are dependent,  $ T_n  \to \infty$ in probability. Such asymptotic theory for $T_n$ was established in the large sample asymptotics where the two dimensions $(p, q)$ are fixed while the sample size $n$ tends to infinity.

When the dimensions $(p,q)$ of the two vectors become large,
 \cite{S13} observed that the above test becomes invalid due to a non negligible bias of the  squared sample distance covariance $\mV_n^2(\bx,\by)$ and then
 proposed a
bias-corrected version $\tilde\mV_n^2(\bx,\by)$ as a substitution.
By this correction, the {\em sample distance
correlation}  $ \tilde{R}_n(\bx,\by) = \tilde\mV_n(\bx,\by) /
[\tilde\mV_n(\bx,\bx)  \tilde\mV_n(\by,\by)]^{1/2}$ was employed for testing the independence hypothesis, whose null distribution was established in a specific asymptotic scheme where $n$ is
kept fixed while $p$ and $q$ both grow to infinity. And this scheme is  referred as fixed-$n$ asymptotic regime in the following.
However, a recent paper \cite{zhu20} reported that 
even the test based on $ \tilde{R}_n(\bx,\by)$ may loss the power for detecting nonlinear correlations
when all the dimensions $(p,q,n)$ grow to infinity.
In particular, they  demonstrated that for high dimensional vectors, their squared sample distance covariance $\tilde\mV_n^2(\bx,\by)$
is asymptotically equivalent
to the summation of their squared component-wise (linear) cross sample
covariances. This implies that  distance covariance 
can only capture linear correlations in such high dimensional regimes.

{
To seek another possibility for detecting non-linear correlations between $\bx$ and $\by$ when all the dimensions $(p,q,n)$ grow to infinity, we propose in this paper a new random matrix model, called {\em distance covariance matrix} (DCM).
Specifically, denote two data matrices  $\X=(\bx_1, \ldots, \bx_n)$ and $\Y=(\by_1, \ldots, \by_n)$, 
the DCM of $\X$ and $\Y$ is defined as 
\begin{align}\label{defsxy}
  \bS_{xy}\triangleq{\bP}_n\D_x{\bf P}_n\D_y{\bP}_n,
\end{align}
where 
\begin{align}\label{Dxy}
  {\D_x}\triangleq\frac{1}{p}{\X'\X}+\frac{1}{pn}\sum_{i=1}^n||\xx_i||^2 {\I_n},\quad {\D_y}\triangleq \frac{1}{q} {\Y'\Y+\frac{1}{qn}\sum_{i=1}^n||\y_i||^2  \I_n},
\end{align}
and
\begin{align}\label{depn}
{\bf P}_n=\I_n-\frac 1n {\bf 1}_n{\bf 1}_n'
\end{align}
 is a projection matrix.
The distance covariance matrix $\bS_{xy}$ is closely connected to the distance covariance $\mV(\bx,\by)$.
As will be discussed in Section \ref{sec:main}, 
a normalized trace of $\bS_{xy}$ is asymptotically equivalent to 
the empirical distance covariance $ \mathcal{V}_n( \bx,\by)$. 
Therefore, we believe that the spectrum of  $\bS_{xy}$ might contain the  information indicating the non-linear dependence nature  between $\bx$ and $\by$. To this end, we 
 investigate the first order asymptotic behaviors of the whole spectrum of the DCM $\bS_{xy}$ under the two-sample  {\em Mar\v{c}enko-Pastur asymptotic regime} 
\begin{align}\label{generalw}
(n, p, q)\to \infty,\quad (c_{n1},c_{n2}):=(p/n, q/n)\to (c_1, c_2) \in(0,\infty)^2.
\end{align} 
Interestingly, we find that instead of the normalized trace of
$\bS_{xy}$, its largest eigenvalues  have the ability to detect
certain non-linear correlation  between the two high dimensional
random vectors $\bx$ and $\by$.
}

The main contributions of the paper are three-folds.
Our first result shows that  the test statistic $T_n$ developed in  \cite{S07} for the independence    hypothesis $H_0$ degenerates to the unit in the
Mar\v{c}enko-Pastur asymptotic regime, which 
extends  a similar
finding in \cite{S13} in their fixed-$n$ asymptotic regime.  Therefore 
the statistic  $T_n$ could not be applied any more for testing the independence    hypothesis $H_0$   in the Mar\v{c}enko-Pastur  asymptotic regime \eqref{generalw}. 

As the second result of the paper, 
we derive a
deterministic limiting distribution $F$  for the eigenvalue distribution of $\bS_{xy}$.
This means in particular that arbitrary eigenvalue statistic of the form
$n^{-1} \sum_i g(\lambda_i)$, where $(\lambda_i)$ denotes the eigenvalues of $\bS_{xy}$, with some smooth function $g$
converges to  $\int g(x)dF(x)$. 
The limiting distribution $F$  is valid  when the vectors
are independent or weakly dependent corresponding  to a
finite-rank perturbation of the independence.
An important property is that this limit
is {\em universal} in the
sense that it does not depend on the details  of the
respective distributions of the  vectors.

%

Third, to  further demonstrate the usefulness of such limiting eigenvalue
distribution, we apply the theory to  the  problem of  detection for certain deviation
from the independence hypothesis by considering a family of
  finite-rank nonlinear dependence alternatives.
{
  We investigate both the global and local spectral behaviors of $\bS_{xy}$.
Globally, because the dependence is of finite rank, the limiting
distribution of the eigenvalues remains the same
as in the independence case,
that is, the universal limit.
However at a local scale,
the  largest eigenvalues of  $\bS_{xy}$  will converge to some limits
outside the support of this universal limit as long as the strength of the dependence
is beyond some critical value.
Moreover,
the locations of these outlying  limits  can be  completely
determined through the model parameters. 
Actually, these results under the finite-rank dependence
  is parallel with what is now known as Baik-Ben-Arous-P\'ech\'e transition in
  random matrix theory,
  see \citet{BBP05},
  \citet{BaikSilv06}
  and \citet{Paul07}.
  In this way, we conclude that the largest eigenvalues of  $\bS_{xy}$
  can be used to detect such dependence structure. 
  In addition, we propose an estimator for the rank of the dependence.
  This estimator is based on the ratios of adjacent largest
  eigenvalues of $\bS_{xy}$. 
  Its performance is assessed through
  simulation experiments.}

Technically,  our theoretical strategy for deriving the universal
limit under independence  is to  derive
a system of equations for  the corresponding  Stieltjes transform
in the Gaussian case first. Indeed when the vectors $\bx$ and $\by$
are Gaussian, the distance covariance matrix $\bS_{xy}$ is
orthogonally invariant; we can thus assume without loss
of generality that the two population covariance matrices are
diagonal,  which  greatly  simplifies the analysis. In
a second step, we obtain an accurate estimate for
the difference between the Stieltjes transforms from Gaussian vectors
and non-Gaussian ones by using a generalization of Lindeberg's substitution
method.
This difference is indeed small enough so
that
the limiting distribution for the global spectrum of $\bS_{xy}$ is
actually universal, regardless of the underlying distributions of the vectors.

%
%
  The rest of the paper is organized as follows.  Section~\ref{sec:main}
  details our model assumptions and  the relation between the  distance covariance matrix  $\bS_{xy}$
and the sample distance covariance $\mV_n(\bx,\by)$.
  Section \ref{sec:lsd}
  establishes the limiting spectral distribution of $\bS_{xy}$ under the
Mar\v{c}enko-Pastur  asymptotic  regime \eqref{generalw}  when $\bx$ and $\by$ are  independent.
  Section \ref{sec:spike} applies this theory to the detection of
  finite-rank nonlinear
  dependence between two
  high-dimensional vectors.
All proofs of our technical results are gathered in  an on-line supplementary file.
  

%
\section{Distance covariance matrix}
\label{sec:main}

Let ${\mathbf M}_p$ be a  $p\times p$ symmetric or Hermitian matrix with eigenvalues $(\lambda_j)_{1\le j\le p}$.
Its spectral distribution  is the probability measure
\[ 
  F^{{\mathbf M}_p}=\frac{1}{p}\sum_{j=1}^p\delta_{\lambda_j}, 
\] 
where  $\delta_b$ denotes  the Dirac mass at $b$. 
For a probability  measure $\mu$ on the real line (equipped with its
  Borel $\sigma$-algebra), its Stieltjes transform
$s_\mu$ is a map from $\mathbb{C}^+$ onto itself,
  \[  s_\mu(z) = \int_{\mathbb{R}} \frac1{x-z} d\mu(x),\qquad z\in \mathbb{C}^+,
  \]
where $\mathbb C^+\triangleq \{z\in \mathbb C: \Im(z)>0\}$.
\medskip

Our asymptotic study of the spectrum of the DCM $\bS_{xy}$ is developed
under the following assumptions.
\begin{description}
\item[Assumption (a)]
  The dimensions $(n,p,q)$ tend to infinity as in \eqref{generalw}.
\item[Assumption (b)]
 The data matrices $\X=(\bx_i)\in\mathbb R^{p\times n}$ and $\Y=(\y_i)\in\mathbb R^{q\times n}$ admit the following independent components model
  $$\X=\Sig_x^{\frac12}\W_{1}\quad \text{and}\quad \Y=\Sig_y^{\frac12}\W_{2},$$
  where $\Sig_x\in\mathbb R^{p\times p}$ and $\Sig_y\in\mathbb R^{q\times q}$ denotes the population covariance matrices of $\bx$ and $\by$, respectively, and  $(\W_1',\W_2')=(w_{ij})$ is an array of i.i.d.\ random variables  satisfying
  $$\E(w_{11})=0,\quad \E(w_{11}^2)=1,\quad \E |w_{11}|^\gamma<\infty,$$
  for some $\gamma\geq 4$.
\item[Assumption (c)]
  The spectral norms of $(\Sig_x, \Sig_y)$ are uniformly bounded
  and their spectral distributions $(H_{xp}, H_{yq})\triangleq (F^{\Sig_x}, F^{\Sig_y})$ converge weakly
  to two probability distributions $(H_x, H_y)$, which are referred as
  {\em population spectral distributions} (PSD).
\end{description}

Our first result concerns the connection between  our distance covariance matrix $\bS_{xy}$ defined in \eqref{defsxy}  and the  sample distance covariance   $\mV_n(\bx,\by)$ defined in \eqref{eq:Vn}.

\begin{theorem}\label{connection}
Suppose that Assumptions (a)-(c) hold with some $\gamma>5$. Then we have
\begin{equation}
  \label{eq:approx}
\mathcal{V}^2_n( \bx,\by) =\frac{1}{2n^2}\sqrt{\frac{pq}{\gamma_x\gamma_y}}\tr \bS_{xy}+o_p(1).
\end{equation}	
\end{theorem}

Theorem \ref{connection} demonstrates that  the squared sample distance covariance $\mathcal{V}^2_n( \bx,\by)$ is asymptotically equal to
the normalized trace of the DCM $\bS_{xy}$.
As a first application of the DCM $\bS_{xy}$,
  we use this approximation to establish below
  the degeneracy of the test statistic $T_n$ given in \eqref{tnsta} for testing the independence hypothesis  \eqref{hoo} under the Mar\v{c}enko-Pastur  asymptotic  framework.

\begin{theorem}\label{th:Tn}
Suppose that Assumptions (a)-(c) hold with some $\gamma>5$.  Then under the null hypothesis $H_0$, we have $\displaystyle T_n\to 1$ in probability.
\end{theorem}

A simple simulation experiment is conducted to exhibit the degeneracy of
$T_n$ for two independent 
standard normal vectors.  The dimension-to-sample size ratios are fixed to be $p/n=q/n=1/2$, the values of $p~(=q)$ range from $50$ to $400$, and the number of independent replications is $1000$.
As shown in Table~\ref{nvs}, with  the growing of $p$, the empirical mean and standard deviation
of $T_n$ converge to   $1$ and $0$, respectively.
  Consequently, the  test established in
  \citet{S07} using the Chi-squares approximation
   will have a much inflated size
  tending to one when the dimensions are indeed large compared to
the  sample size. 
\begin{table}[H]
 \setlength\tabcolsep{6pt}
 \begin{center}
 \caption{Empirical mean and standard deviation of the test
      statistic $T_n$ from 1000 independent replications with $p/n=q/n=1/2$ and $p\in\{50,100,200,400\}$. Independent standard normal vectors are used for  $\bx$ and $\by$.\label{nvs}}
    \begin{tabular}{ccccccccccc}
       \hline
      \multicolumn{2}{c}{$p=50$}&\multicolumn{2}{c}{$p=100$} &\multicolumn{2}{c}{$p=200$}&\multicolumn{2}{c}{$p=400$}\\
      \hline
      mean&sd&mean&sd&mean&sd&mean&sd\\
      1.0104&0.0075&1.0048&0.0036&1.0026&0.0018&1.0013&0.0009\\
      \hline
    \end{tabular}
 \end{center}
\end{table}


%
\section{Limiting spectral distribution of $\bS_{xy}$ when $\xx$ and $\y$ are independent}
\label{sec:lsd}

This section presents the first order convergence of the empirical spectral distribution $F^{\bS_{xy}}$ of the DCM $\bS_{xy}$  when $\xx$ and $\y$ are independent.
\begin{theorem}\label{lsd}
  Suppose that Assumptions (a)-(c) hold. Then, almost surely, the
  empirical spectral distribution $F^{\bS_{xy}}$ converges weakly to a
  {limiting spectral distribution (LSD)} $F$ whose Stieltjes transform
  $s=s(z)$ is a solution to the following system of equations:
  \begin{align}\label{lsd-sys}
	\begin{cases}
      \displaystyle	s=\frac{wm-1}{z},\\
      \displaystyle	w=\int ts+\frac{ts}{1+tc_1^{-1}sm}dH_x(t),\\
      \displaystyle	m=\int t+\frac{t}{1+tc_2^{-1}w}dH_y(t),
	\end{cases}
  \end{align}
  where $w=w(z)$ and $m=m(z)$ are two auxiliary analytic
  functions. The solution is also unique on  the set
  \begin{align}\label{set}
	\{s(z): s(z)\in \mathbb C^+, w(z)\in \mathbb C^+, m(z)\in \mathbb C^-, z\in \mathbb C^+\}.
  \end{align}
\end{theorem}

\begin{remark}
 The two auxiliary functions $w(z)$ and $m(z)$ are respectively the
 limits of $w_n(z)$ and $m_n(z)$ defined in \eqref{wmn} of the supplementary document.
 Their
 construction accounts for the signs of their imaginary parts as in
  \eqref{set}.
  \end{remark}

 \begin{figure}[h]
\centering
    \includegraphics[width=2.5in]{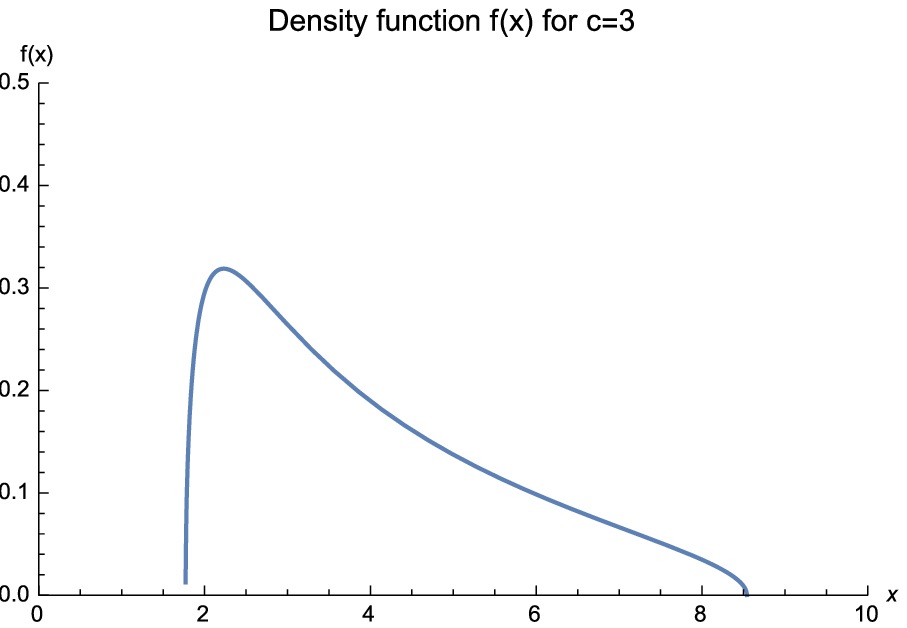}
    \includegraphics[width=2.5in]{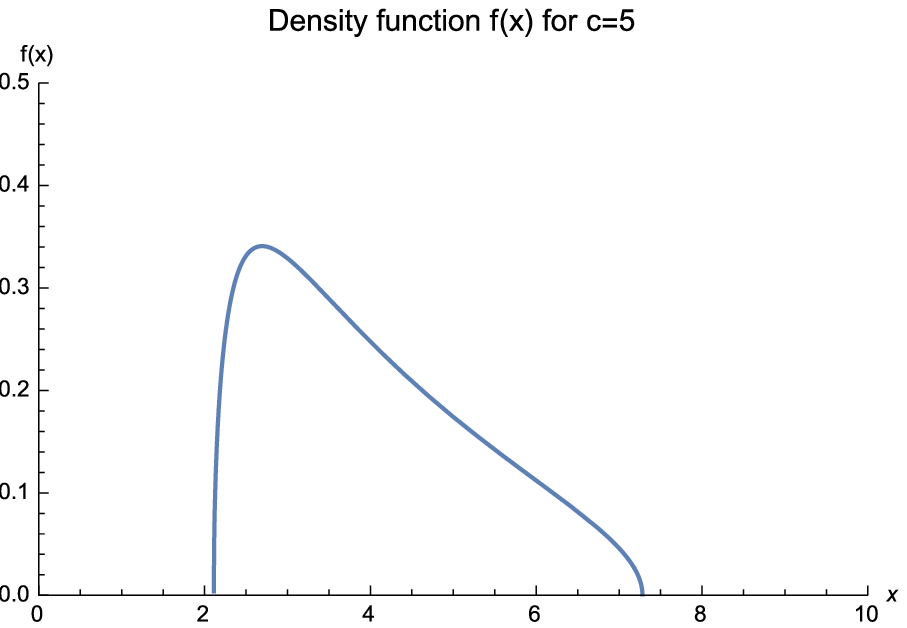}
  \caption{Density curves of LSDs for $c=3$ (left) and $c=5$ (right). The PSDs are $H_x=H_y=\delta_1$.
  }
  \label{lsd35}
\end{figure}

\red{
Next, we provide an illustration of how to calculate the LSD $F$ through the system of equations \eqref{lsd-sys}.
Considering the case where the two populations $\bx$ and $\y$ are of the same dimension and both have the identity covariance matrices, we thus have  
\begin{align}\label{example}
c_1=c_2=c\quad\text{and}\quad H_x=H_y=\delta_1.
\end{align}
For this case, 
a closed-form solution to the system \eqref{lsd-sys}  does exist, that is,
the Stieltjes transform $s=s(z)$ of the LSD $F$ satisfies the following
\begin{align}\label{mp0}
c^2 - s + 2 c s - 4 c^2 s + s^2 + c^2 s z - 2 s^2 z + 2 c s^2 z + 
 s^3 z - s^3 z^2=0.
\end{align}
Substituting $z=x+{\rm i}v$ and $s=s_u+{\rm i}s_v$ into \eqref{mp0} and then letting $v\downarrow0$, we get the following system of equations by separating the real and imaginary parts on the left hand side of \eqref{mp0},
\begin{align}
\left\{
\begin{array}{l}
\displaystyle
s_v^2=\frac{c^2 - s_u + 2 c s_u - 4 c^2 s_u + c^2 x s_u + s_u^2 - 2 x s_u^2 + 
 2 c x s_u^2 + x s_u^3 - x^2 s_u^3}{1 - 2 x + 2 c x + 3 x s_u - 
 3 x^2 s_u},\vspace{1ex}\\\displaystyle
 s_v^2=\frac{1 - 2 c + 4 c^2 - c^2 x - 2 s_u + 4 x s_u - 4 c x s_u - 3 x s_u^2 + 
 3 x^2 s_u^2}{-x + x^2}. \label{sv1}
 \end{array}
 \right.
\end{align}
Cancelling the variable $s_u$ from \eqref{sv1}, one may get three solutions for $s_v^2$ as a function of $x$. These three functions indeed have closed forms, but are lengthy and we omit their explicit expressions here. Then  for each real value of $x$, only one solution of $s_v^2$ is real and nonnegative, which corresponds to the density function $f(x)$ of the LSD $F$, i.e. $f(x)=\sqrt{s_v^2}/\pi$. Using  this approach, we plot in the following Figure \ref{lsd35} two LSDs for such particular setting \eqref{example} corresponding to $c=3$ and $c=5$.}
 However, generally when there is no closed-form solution  for  \eqref{lsd-sys}, we 
   rely on numerical approximations for the limiting Stieltjes
   transform $s(z)$ and the underlying limiting density
   function. These methods are used in the illustration below and also
   in the simulation experiments in Section~\ref{sec:spike}.

Some numerical illustrations of  Theorem \ref{lsd} are conducted under two models:
\begin{itemize}
\item []Model 1: $H_x=H_y=\delta_1$, $c_1=c_2=1$, $z_{11}\sim N(0,1)$;
\item []Model 2: $H_{x}=0.5\delta_{0.5}+0.5\delta_1$,
  $H_{y}=0.5\delta_{0.25}+0.5\delta_{0.75}$, $c_1=2, c_2=1$, and
  $z_{11}\sim (\chi^2_v-v)/\sqrt{2v}$, a standardized chi-squared
  distribution  with degree of freedom $v=2$.
\end{itemize}
The PSDs in the  first model are simple point masses and  the system
\eqref{lsd-sys} defining the LSD  simplifies  to a single equation
$(z^2-z)s^3-s^2+(3-z)s-1=0$ (letting $c=1$ in \eqref{mp0}).
The second model is a bit more elaborated where  the PSDs are mixtures
of  two point  masses and the innovations  $z_{ij}$'s  are  Chi-square
distributed with heavy tails.

%

To exhibit the LSDs defined by Models 1 and 2, we simply approximate their density functions by
$\hat f(x) =\Im s(x+{\rm i}/ 10^{4})/\pi,  x\in \mathbb R.$
This approximation is justified by  the inversion formula of Stieltjes transforms, i.e.
$f(x)=\lim_{\varepsilon\to 0^+} \Im s(x+{\rm i}\varepsilon)/\pi,$
provided  the limit exists, see Theorem B.10 in \cite{BSbook}.
Obviously, our approximation takes $\varepsilon=10^{-4}$ which is
small enough for the illustration here. Next, for any given $z=x+{\rm i}/ 10^{4}$, we numerically solve the system of equations in \eqref{lsd-sys} and select the unique solution $(s(z),w(z), m(z))$ satisfying \eqref{set}, which is done automatically in Mathematica software. Finally, taking the imaginary part of $s(z)/\pi$ gives $\hat f(x)$.

 \begin{figure}
	\includegraphics[width=2.7in]{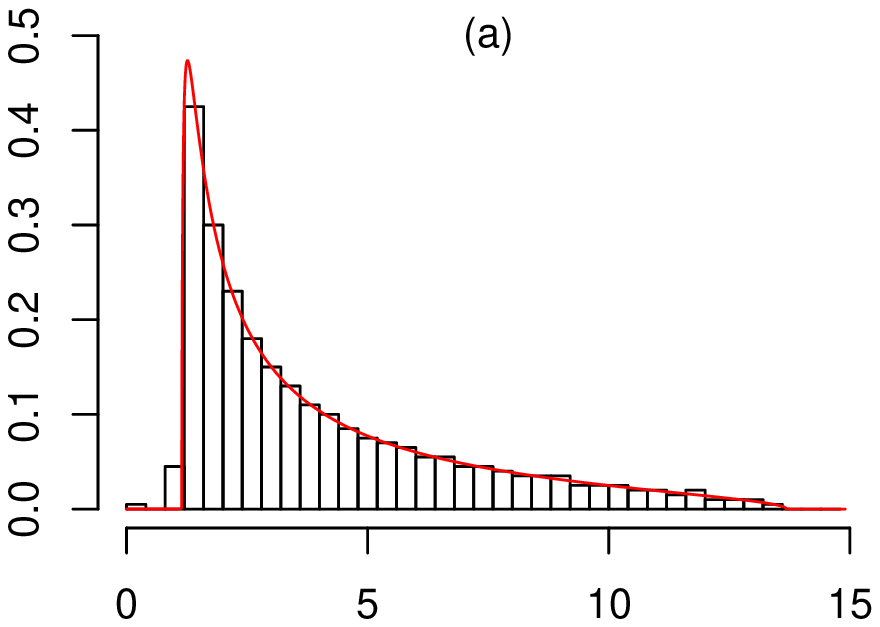}
	\includegraphics[width=2.7in]{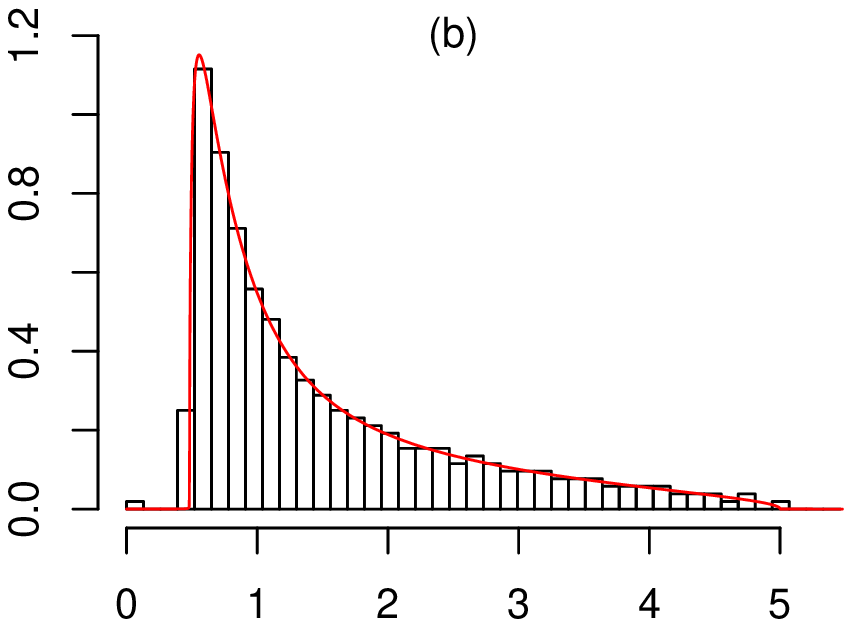}
   \caption{Histogram of eigenvalues of the matrix $\bS_{xy}$ under Model 1 (left panel) with dimensions $p=q=n=500$ and Model 2 (right panel) with dimensions $p=2q=2n=800$. The fitted curves (red colored) are the corresponding densities of the LSDs.}
   \label{dens}
\end{figure}

In this simulation experiment, the empirical PSDs are chosen as their
limiting ones and the dimensions are $(p,q,n)=(500,500,500)$
for Model 1 and $(p,q,n)=(800,400,400)$ for Model 2. All eigenvalues
are collected from 100 independent replications. The averaged
histograms of the eigenvalues of $\bS_{xy}$ from those replications are
depicted in Figure \ref{dens}.  It shows that these empirical
distributions  match well    their limiting density curves predicted in Theorem \ref{lsd}.

\section{Application to the detection of  dependence between two high dimensional vectors}\label{sec:spike}

Theorem~\ref{lsd} determines  a universal limit for the  bulk
spectrum  of  the DCM
when the two sets of samples are
independent. A natural question arises:  how this bulk limit
will evolve when they
become dependent?
Apparently if their inherent dependence  is very strong, the spectral
limit of  the DCM will be totally different from the universal limit
in Theorem~\ref{lsd}.
  Here we  choose to study   a special type of  weak dependence which is finite rank dependence. Such concept is parallel to the idea of finite-rank perturbation or  spiked population
  models in high-dimensional statistics
  which are  widely studied
  in
  connection with high-dimensional PCA, factor modeling and the signal
  detection problem \citep{JohnstonePaul18}.
A striking finding from the work here is that such finite-rank
nonlinear dependence
can be detected using the largest eigenvalues of the DCM
while existing methods based on sample covariance, sample correlation
or sample canonical correlations will fail.

\subsection{Extreme eigenvalues of distance covariance matrix    under  finite-rank dependence}
Precisely, we consider two dependent populations  $\xx  \in
\mathbb R^p$ and $\z\in \mathbb R^q$ defined as follows:

\begin{itemize}
\item[(i)] For a fixed $m\in\mathbb N$, let $(\bu_k)_{1\leq k\leq
    m}$ and $(\bv_k)_{1\leq k\leq m}$ be two independent sequences of
  i.i.d.\ vectors distributed uniformly  on the unit  spheres in $ \mathbb R^q$ and
  $\mathbb R^p$, respectively.
\item[(ii)] Given the sequences $(\bu_k)$ and $(\bv_k)$,  the population $\z$ is defined as
\begin{align}\label{zx}
  \z=\varepsilon \left(\sum_{k=1}^m \theta_k \bu_k \bv'_k\right)\xx+\y,
\end{align}
\end{itemize}
where
\begin{itemize}
\item[(1)]  $\xx  \in \mathbb R^p$ and $\y\in \mathbb R^q$ satisfy Assumptions (b) and (c);
\item[(2)]  \red{$\varepsilon$ is a standardized random variable with finite fourth moment.}
\item[(3)]    $0<\theta_m<\cdots<\theta_1<\infty$ are $m$ constants reflecting the  strengths of dependence  between $\xx$ and $\z$.
\end{itemize}

\begin{remark}
  The pair of random vectors $(\xx,\z)$ in \eqref{zx} are nonlinearly
  dependent, that
  is, they are uncorrelated but dependent.  To see this, consider a particular case such that $\varepsilon$ is a random sign  taking  values $1$ or $-1$ with equal probability.  Then  it is easy to see
  that the random sign $\varepsilon$ put on the vector $\bx$ implies
  the uncorrelation between the vectors. To establish their
  dependence, simple algebra shows that
  \begin{align*}
    \E \big( \|\xx\|^2\big)\E \big( \|\z\|^2\big)
    & =  \frac1p \|\theta\|^2  \E^2 \big( \|\xx\|^2\big) + \E\|\xx\|^2\E\|\y\|^2,\\
    \E \big( \|\xx\|^2\|\z\|^2 \big)
    & =  \frac1p \|\theta\|^2  \E \|\xx\|^4 + \E\|\xx\|^2\E\|\y\|^2.
  \end{align*}
  Here, $\|\theta\|^2 =\theta_1^2+\cdots+\theta_m^2$.
  Unless $\xx$ is a constant vector, $\E \|\xx\|^4>\E^2 \big(
  \|\xx\|^2\big)$ and thus the vectors $\xx$ and $\z$ are dependent.
\end{remark}

Suppose we have  an i.i.d.\  sample  $(\xx_1,\z_1),\ldots,(\xx_n,\z_n)$ from the population $(\xx,\z)\in\mathbb R^p\times \mathbb R^q$  defined in \eqref{zx}. Denote by $\X=(\xx_1,\ldots,\xx_n)$ and $\Z=(\z_1,\ldots,\z_n)$ the two data matrices with sizes $p\times n$ and $q\times n$, respectively.
Similar to the matrices in \eqref{Dxy} we define two matrices $\D_x$ and $\D_z$ as
\begin{align*}
\D_x=\frac{1}{p}{\X'\X}+\kappa_{x}{\bf I}_n\quad  \text{and} \quad\D_z=\frac{1}{q}{\Z'\Z}+\kappa_{z}{\bf I}_n,
\end{align*}
where $$\kappa_{x}\triangleq(pn)^{-1}\sum_{i=1}^n||\bx_i||^2\quad \text{and}\quad\kappa_{z}\triangleq(qn)^{-1}\sum_{i=1}^n||\z_i||^2.$$ The corresponding DCM is  written as $$\bS_{xz}\triangleq\bP_n \D_x \bP_n\D_z\bP_n.$$
We will study  the  spectral properties of $\bS_{xz}$ for the  dependent
pair $(\xx,\z)$ defined in \eqref{zx}.
  First of all, because the rank of perturbation is finite,
  it is shown that
  the limiting spectral distribution of $\bS_{xz}$ remains the
  same as if the two populations are independent.

\begin{theorem}\label{lsd-spike}
 Suppose that Assumptions (a)-(c) hold for model \eqref{zx}. The limiting spectral distribution of $\bS_{xz}$ is given by the same $F$ defined in Theorem~\ref{lsd}.
\end{theorem}
According to Theorem \ref{lsd-spike}, the global behavior of the eigenvalues of the DCM $\bS_{xz}$ will not be
  useful for distinguishing such weak dependence from the
  independence scenario.
  In the following, we turn to study the top eigenvalues of the $\bS_{xz}$ and
  show that the weak dependence structure is encoded in these top eigenvalues.
   Detection of such weak dependence thus  becomes possible using
  these top eigenvalues.
Before that, we introduce some notations  that will be used for stating our result.
We denote $$\lambda_+=\limsup_{n \to \infty}\|\bS_{xy}\|,$$ which is finite.
On $(\lambda_+,\infty)$, define the function
  \begin{align}\label{function-g}
    g(\lambda)=-\int tdH_x(t)\int \frac{w(\lambda)}{c_2+t w(\lambda)}dH_y(t),\quad \lambda>\lambda_+,
  \end{align}
  where  $w(z)$ is given in \eqref{lsd}.
  It's easy to verify that
  $g(\lambda)>0$,  $g'(\lambda)<0$ and  $\lim_{\lambda\to+\infty}g(\lambda)= 0.$
Next, define
  \begin{align}\label{theta0}
    \theta_0 := \lim_{\lambda\,\downarrow\, \lambda_+} [g(\lambda)]^{-\frac 12}.
  \end{align}
  Therefore $g$ is a one-to-one,  strictly decreasing and nonnegative
  function from
  $(\lambda_+,\infty)$ to $(1/\theta_0^2,0)$.

  \begin{theorem}\label{th-spike}
   Suppose that  Assumptions (a)-(c) hold  for model \eqref{zx} and
    for some  $k \in \{1,\ldots, m\}$,
    $\theta_k>\theta_0$.  Then the $k$-{\em th} largest
    eigenvalue $\lambda_{n, k}$ of the DCM $\bS_{xz}$
    converges almost surely to a limit
    \begin{align}\label{outlier}
      \lambda_k = g^{-1}\big( 1/{\theta_k^2}\big) > \lambda_+,
    \end{align}
where $g^{-1}$ denotes the functional inverse of $g$.
  \end{theorem}
  
\red{
\begin{remark}
Generally,  the function $g$ as well as the critical value $\theta_0$ have no analytic formulas and  both can be found numerically for any given model settings. However, in certain particular case, for example 
 the setting considered in \eqref{example}, the function $g(\lambda)$ given in \eqref{function-g} is a solution to
\begin{align}\label{gf}
c g^3(\lambda)  + (1 + 4 c) g^2(\lambda) +  g(\lambda) (3 + 4c - c\lambda)+2=0.
\end{align}
In fact there are three solutions to \eqref{gf} which all have explicit but lengthy expressions. 
One may choose the one that monotonically decreases to zero as $\lambda$ tends to infinity, which is our target function $g(\lambda)$. Then the critical value $\theta_0=[g(\lambda_+)]^{-1/2}$ can be obtained accordingly. Note that the right edge $\lambda_+$ of the LSD $F$ can  be   theoretically derived by setting the  density function $f(x)$ to be zero. 
 As an illustration, we exhibit the relation between the value $\theta_0$ and the ratio $c$  for the case \eqref{example}  in Figure \ref{theta00}.
\begin{figure}
\centering
    \includegraphics[width=2.8in]{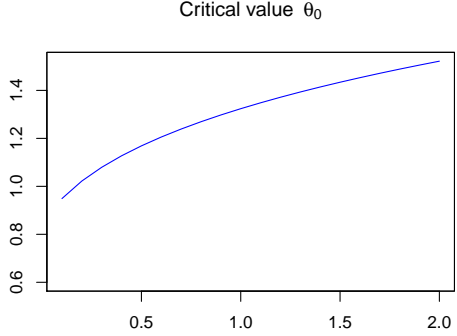}
    \setlength{\abovecaptionskip}{-0.5cm} 
  \caption{Critical value $\theta_0$ for $c_1=c_2=c$ and $H_x=H_y=\delta_1$.
  }
  \label{theta00}
\end{figure}
\end{remark}
}

  The limit $\lambda_k$ in \eqref{outlier} is  outside the support of the LSD $F$.
  A technical point here is that
  Theorem~\ref{th-spike} does not tell
  what happens to $\lambda_{n,  k}$  if
  $\theta_k\le \theta_0$.
  By assuming the convergence of the largest eigenvalue of the
  base component $\bS_{xy}$ to the right edge point of the
  LSD, we can establish  the following {\em exact phase transition}
  for the top eigenvalues  $\lambda_{n, k}~(1\le k\le m)$.

\begin{corollary}\label{co-spike}
  In addition to Assumptions (a)-(c) for model \eqref{zx}, suppose that the largest
  eigenvalue of the DCM $\bS_{xy}$ converges to $\lambda_+$, which is the right edge point
  of the LSD $F$.
  Then for  $k=1,\ldots,m$,
$$
\lambda_{n,k}\xrightarrow{a.s.}
\begin{cases}
\lambda_k& \text{if}~ \theta_k>\theta_0,\\
\lambda_+& \text{if}\ \theta_k\leq\theta_0,
\end{cases}
$$
where $\theta_0$ and $\lambda_k$ are  given in \eqref{theta0} and
\eqref{outlier}, respectively.
\end{corollary}
Corollary \ref{co-spike} follows directly from the proof of Theorem
\ref{th-spike} and the classic interlacing theorem.
It implies the value $\theta_0$ is the exact critical value for
  the phase transition of the top eigenvalues of the DCM $\bS_{xz}$.
  Note that the convergence of the largest eigenvalue of the (null)
  DCM  $\bS_{xy}$ to $\lambda_+$ is needed and assumed here to ensure the
  convergence of those sub-critical  spike eigenvalues, i.e.\
  $\theta_k\le \theta_0$, to the same right edge point
  $\lambda_+$. On the other hand, very likely this largest eigenvalue
  does converge. However,
  proof for such convergence of
  the largest eigenvalue
  would be  lengthy and technical, and  we leave it for future investigation.

\subsection{Monte Carlo experiments}\label{simu}
This section examines  finite sample properties  of the outlier eigenvalues of $\bS_{x z}$.
To simplify the exposition, we consider only the rank-one
  situation  ($m=1$) in this section. Higher dependence ranks with $m>1$ will be
  discussed in Section~\ref{estimating-m}.
Three models are taken into consideration under normal populations:\\
{\bf Model 4}:
  $
  H_x=H_y=\delta_1,\ c_1=c_2=2;\\
  $
{\bf Model 5}:
  $
  H_x=H_y=\delta_1,\ c_1=0.1,\ c_2=0.2;\\
  $
{\bf Model 6}:
  $
  H_x=0.5\delta_{0.5}+0.5\delta_{1},\ H_y=0.5\delta_{1}+0.5\delta_{1.5},\ c_1=1, c_2=2.\\
  $
  Models 4 and 5 are both standard normal population,  with different dimension-to-sample size ratios. Model 6 is more general by employing two discrete PSDs. All statistics are calculated using 1000 independent replications.

We begin with  the convergence of the largest eigenvalue of $\bS_{xz}$
under Model 4. Theoretically, the largest eigenvalue will become an
outlier when $\theta>{\  \theta_0=}1.52$ (see Figure \ref{theta00} for the critical value).  The parameter $\theta$ is
thus set to be $\theta=0,1,2,3$. The sample size $n$ ranges from 100
to 1600. Empirical mean and standard deviation of the largest
eigenvalue are collected in Table \ref{table11}. It shows that, for
$\theta=0$ and 1 (second to fifth columns), the largest eigenvalue
increases with decreasing standard error as $n$ grows and is close to
$\lambda_+=9.95$, the right edge point of  $F$. When  $\theta={
  2}$ and { 3}
(last four columns), the largest eigenvalue converges to its
theoretical limit $\lambda=10.6875$ for $\theta=2$ and
$\lambda=15.0123$ for $\theta=3$.  These results fully coincide with  the conclusions of Theorem \ref{th-spike}.

\begin{table}[H]
    \caption{Empirical mean and standard deviation of the largest eigenvalue under Model 4. The setting is
      $c_{n1}=c_{n2}=2$ with varying  $n$ and 
      1000 independent replications. The right edge point of the
      LSD is $\lambda_{+}=9.95$.\label{table11}}
    \centering
    \resizebox{\textwidth}{!}{
    \begin{tabular}{cccccccccccc}
          \hline
      & \multicolumn{2}{c}{$\theta=0$}&\multicolumn{2}{c}{$\theta=1$} &\multicolumn{2}{c}{$(\theta,\lambda)=(2,10.6875)$}&\multicolumn{2}{c}{$(\theta,\lambda)=(3,15.0123)$}\\
      \hline
      $n$&mean&sd&mean&sd&mean&sd&mean&sd\\
      \cline{2-3} \cline{4-5}\cline{6-7}\cline{8-9}\cline{10-11}
      $100$ &9.5732&0.3126&9.6443&0.3419&10.7285&0.7055&15.1056&1.5770\\
      $200$ &9.7247&0.1972&9.7486&0.2013&10.7219&0.5048&15.0821&1.1099\\
      $400$ &9.8094&0.1302&9.8209&0.1239&10.7114&0.3500&15.0446&0.7458\\
      $800$ &9.8587&0.0769&9.8729&0.0796&10.7079&0.2531&14.9985&0.5505\\
      $1600$&9.8950&0.0479&9.8966&0.0502&10.6985&0.1745&15.0249&0.3794\\
      \hline
    \end{tabular}
    }
    \end{table}

\begin{figure}
  \phantom{abcdefg}
  \vspace{-0.5cm}
  \centering
    \includegraphics[width=0.45\linewidth]{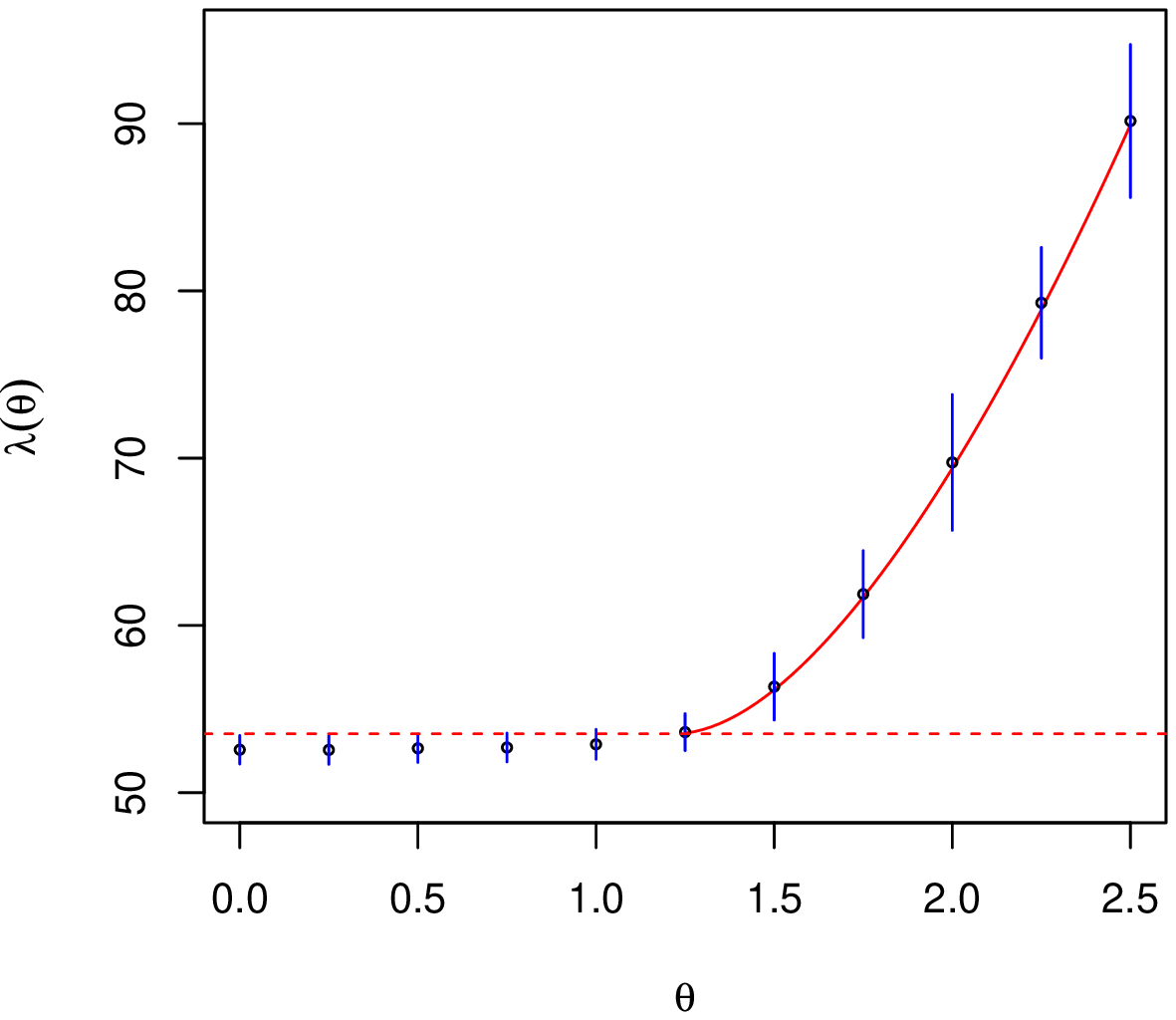}
    \includegraphics[width=0.45\linewidth]{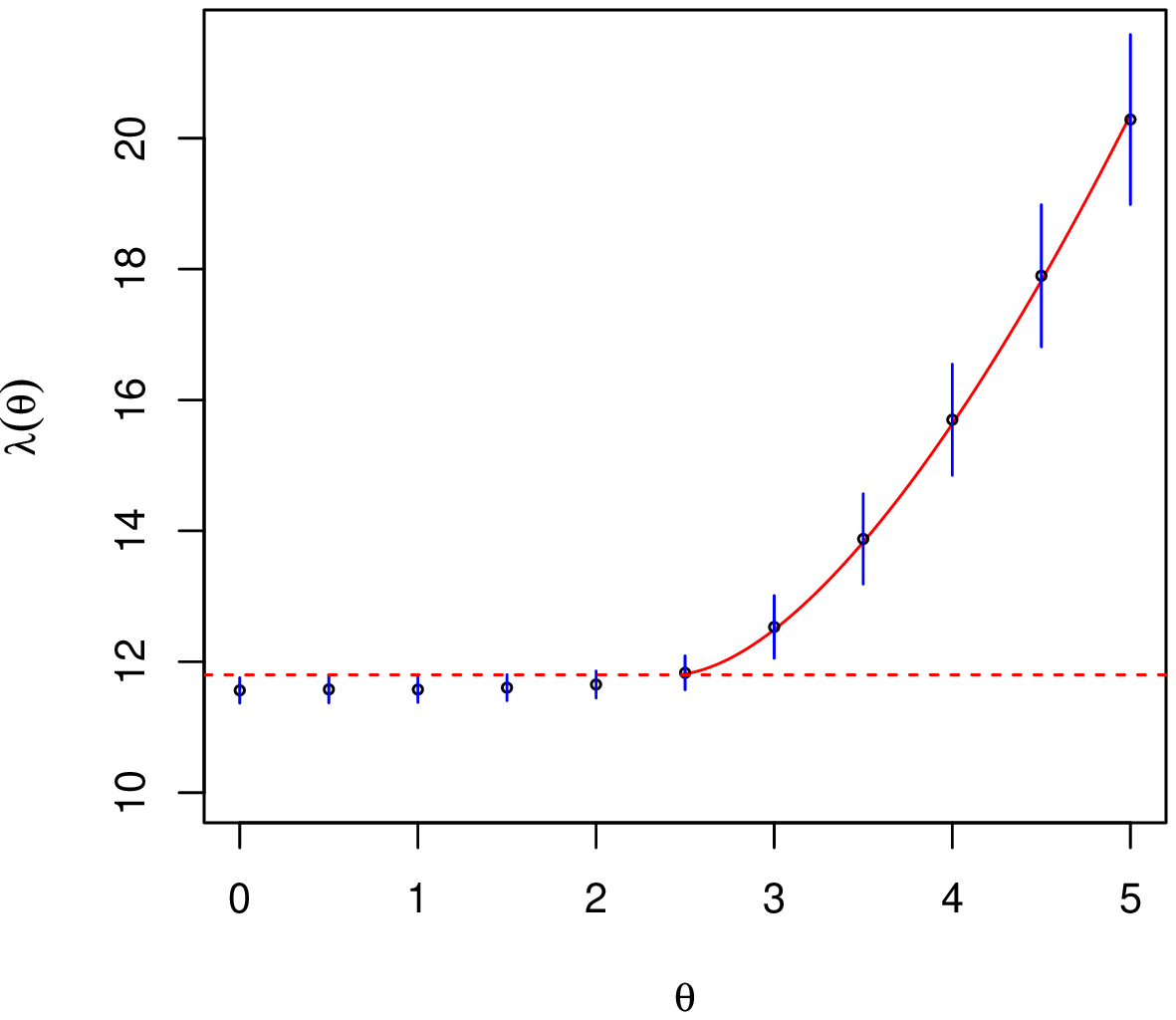}
    \setlength{\abovecaptionskip}{-0.4cm} 
  \caption{The average of the largest eigenvalue under Model 5 (left
    panel) and Model 6 (right panel) from 1000 independent
    replications with $\pm 1$ standard deviations (blue bars). The
    solid red line is the limiting curve of the function
    $\lambda(\theta)$ and the dashed red line represents the right
    boundary of LSD's support.
  }
  \label{spike-1}
\end{figure}

Next we study the evolution of the outlier limit  $\lambda(\theta)$ in
function of the dependence strength $\theta$.  Models 5 and 6 are
considered with  the dimensions  fixed at $(p,q,n)=(200, 400, 2000)$ for
Model 5 and at $(p,q,n)=(800, 800, 400)$ for Model 6. The parameter
$\theta$ ranges from 0 to 2.5 for Model 5 and from 0 to 5 for Model
6. Figure \ref{spike-1} displays the average of the largest eigenvalue
with $\pm1$ standard deviations (vertical bars). The dashed red lines
mark the right boundary of $F$ and the solid red lines are the
theoretical curves of $\lambda=\lambda(\theta)$. Both the two graphs
in Figure \ref{spike-1} exhibit a common trend that the largest
eigenvalue will depart from the bulk when $\theta$ crosses a  critical value and goes up with increasing standard deviation.

\begin{figure}
    \begin{minipage}[t]{0.33\linewidth}
	\includegraphics[width=1.8in]{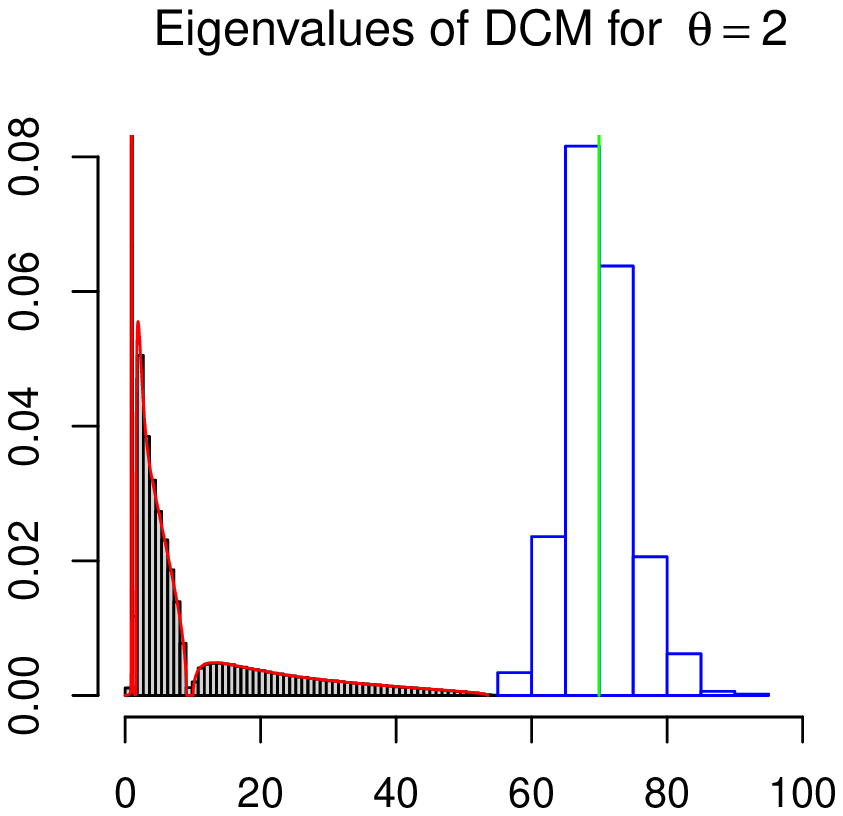}
  \end{minipage}%
  \begin{minipage}[t]{0.33\linewidth}
	\includegraphics[width=1.8in]{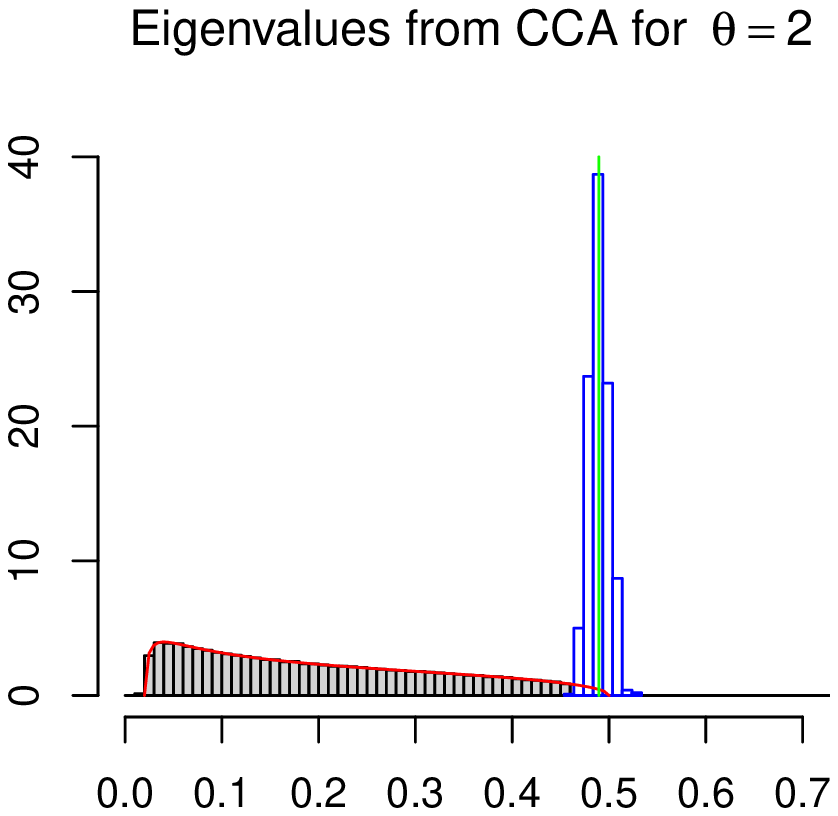}
  \end{minipage}
   \begin{minipage}[t]{0.33\linewidth}
	\includegraphics[width=1.8in]{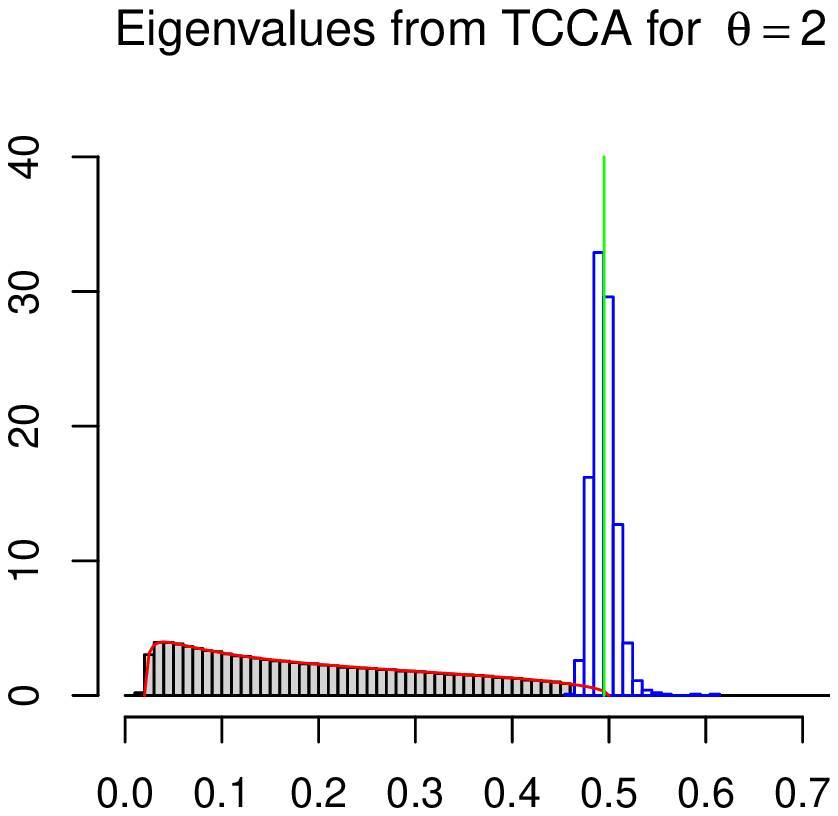}
  \end{minipage}
      \begin{minipage}[t]{0.33\linewidth}
	\includegraphics[width=1.8in]{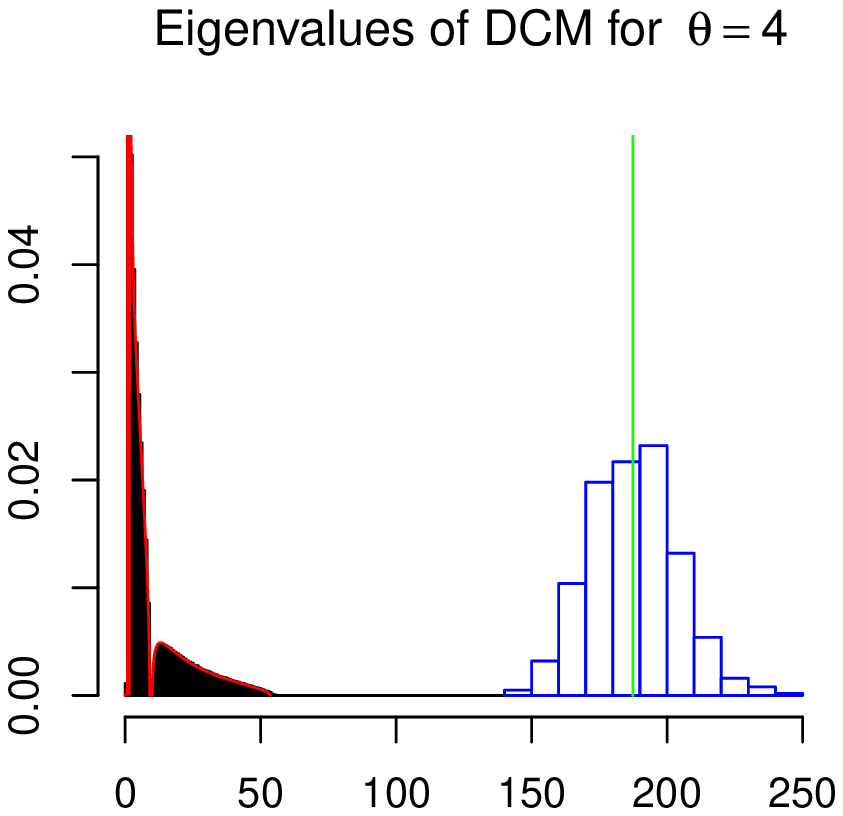}
  \end{minipage}%
  \begin{minipage}[t]{0.33\linewidth}
	\includegraphics[width=1.8in]{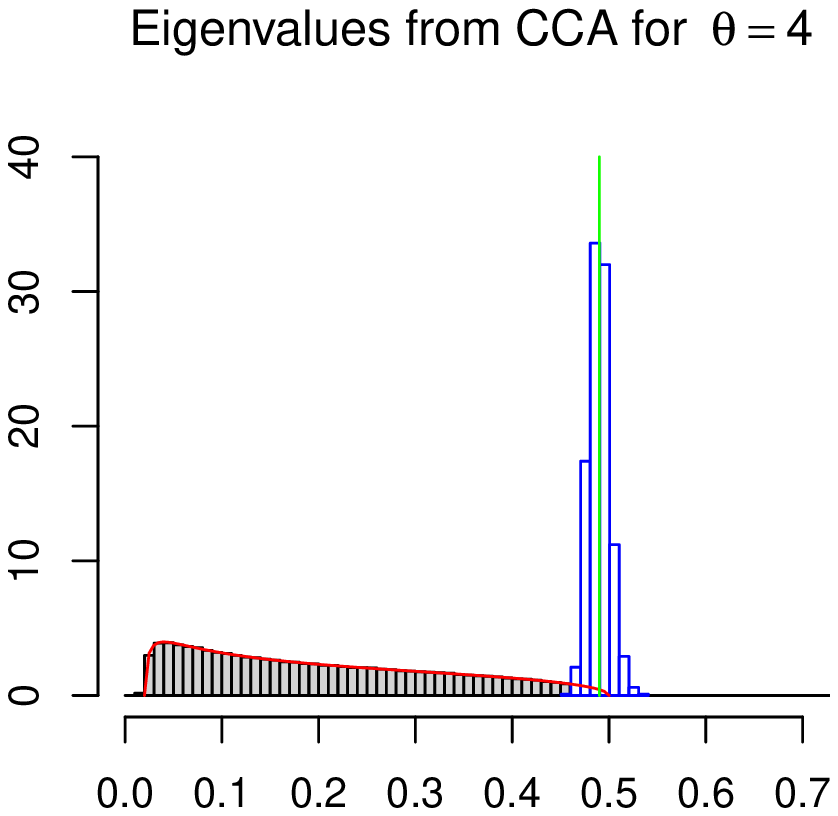}
  \end{minipage}
   \begin{minipage}[t]{0.33\linewidth}
	\includegraphics[width=1.8in]{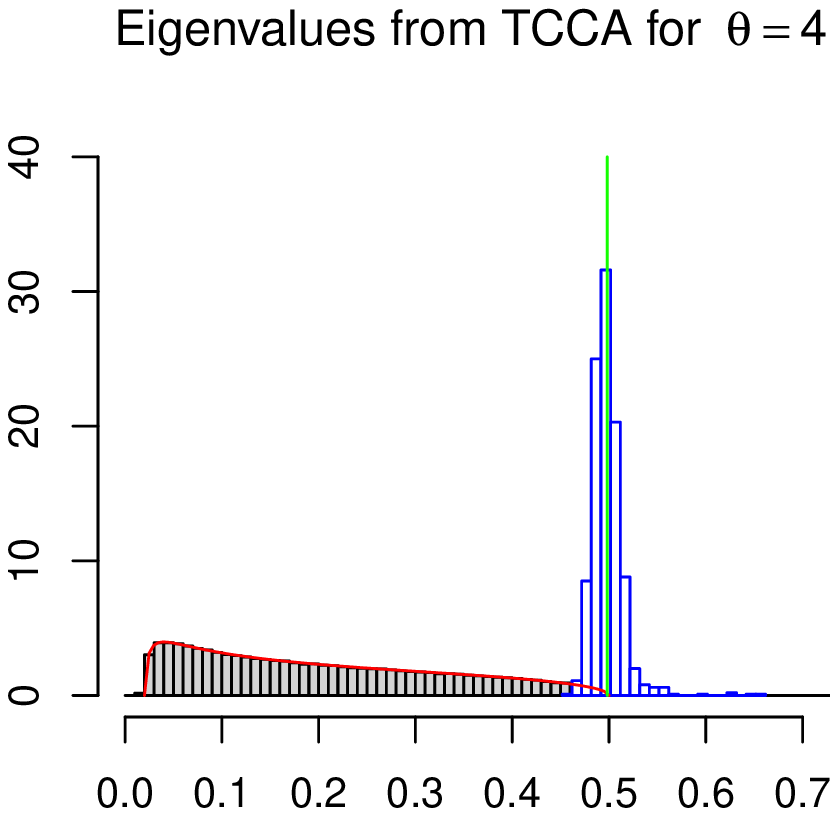}
  \end{minipage}
    \begin{minipage}[t]{0.33\linewidth}
	\includegraphics[width=1.8in]{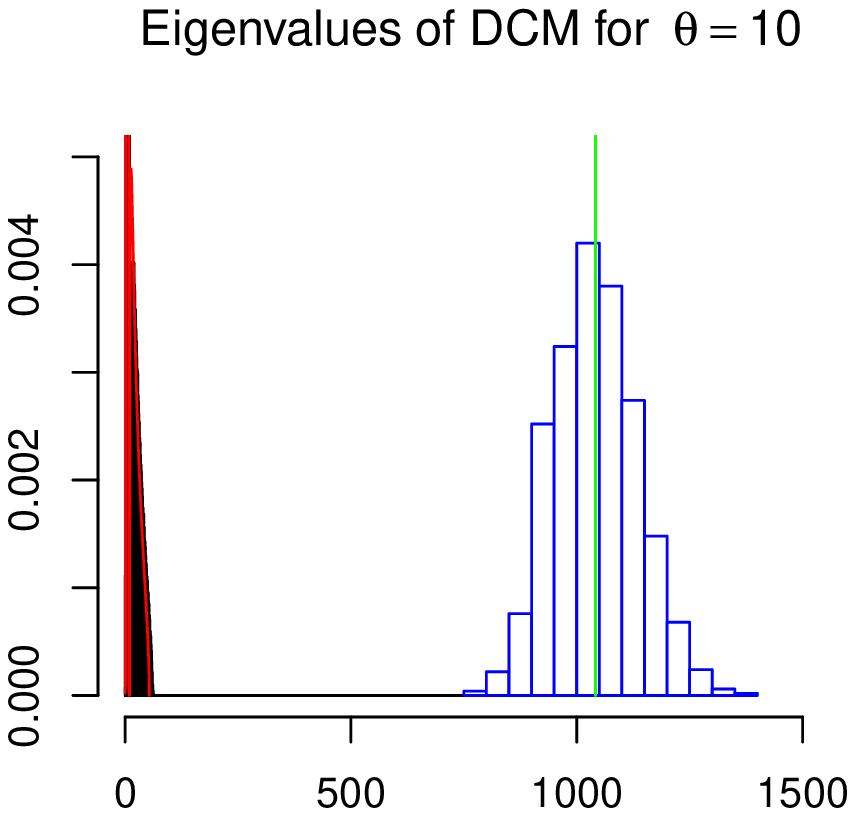}
  \end{minipage}%
  \begin{minipage}[t]{0.33\linewidth}
	\includegraphics[width=1.8in]{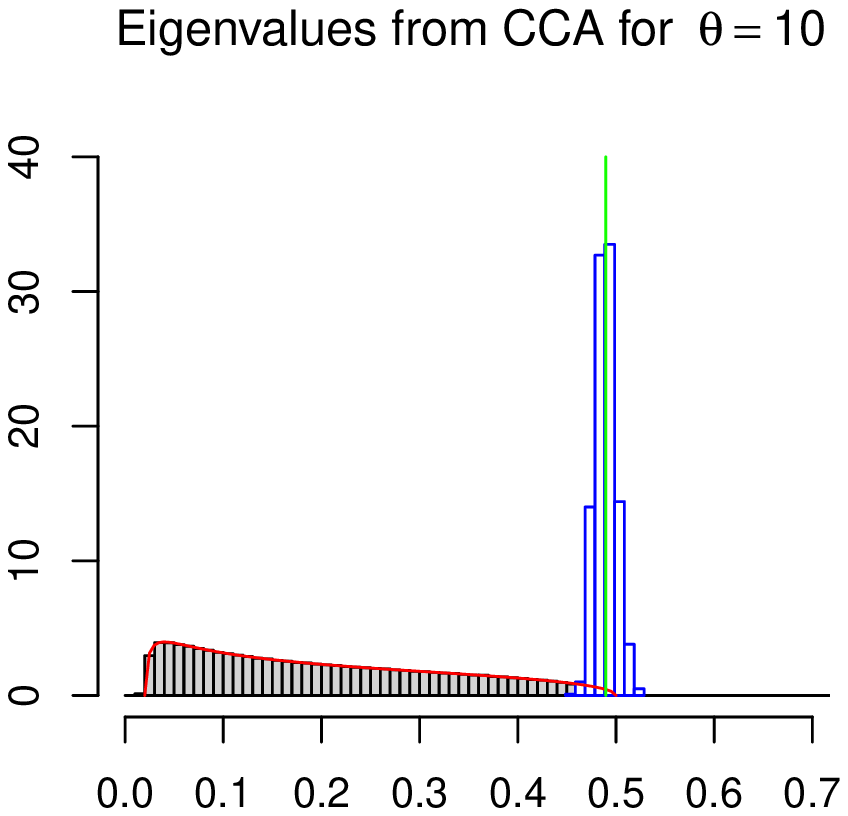}
  \end{minipage}
   \begin{minipage}[t]{0.33\linewidth}
	\includegraphics[width=1.8in]{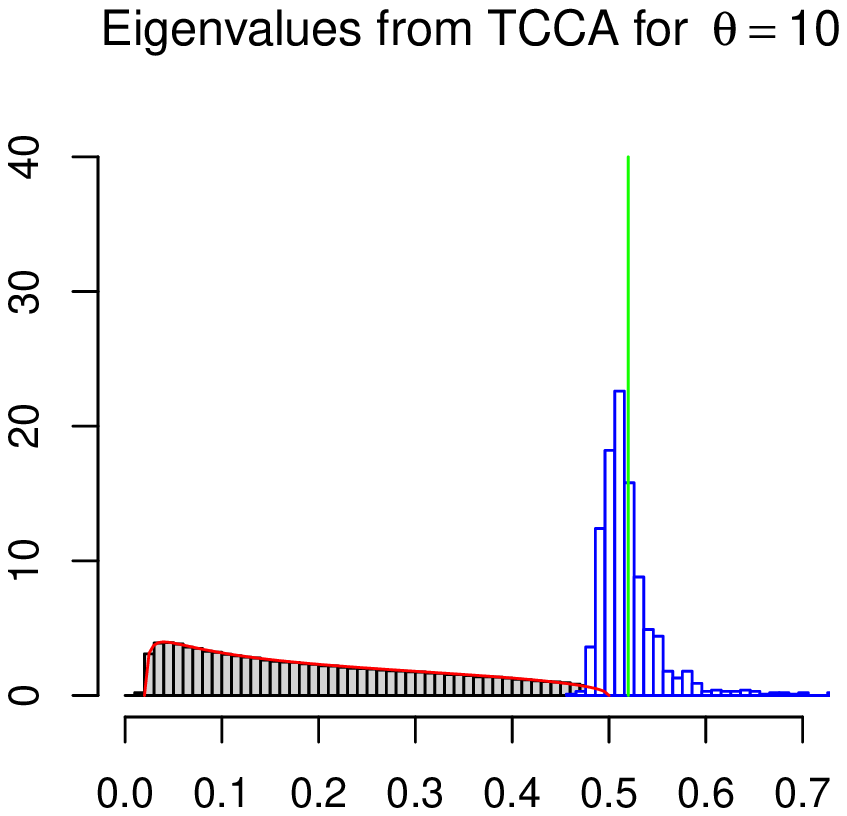}
  \end{minipage}
  \setlength{\abovecaptionskip}{-1.8cm} 
  \caption{Histograms of bulk eigenvalues (black strips) and the largest eigenvalue (blue strips) from 1000 independent replications under Model 5 with $\theta=2,4,$ and $10$. The red curves are LSD densities and the green vertical lines locate at the averages of the largest eigenvalues.  Plots on the left panel are based on the DCM $\bS_{xz}$ while those on the middle and left panels are based on CCA and TCCA, respectively. The dimensions are $(p,q,n)=(100,200,1000)$.}    \label{spike-cca}
\end{figure}

\begin{figure}
	\begin{minipage}[t]{0.33\linewidth}
		\includegraphics[width=2in]{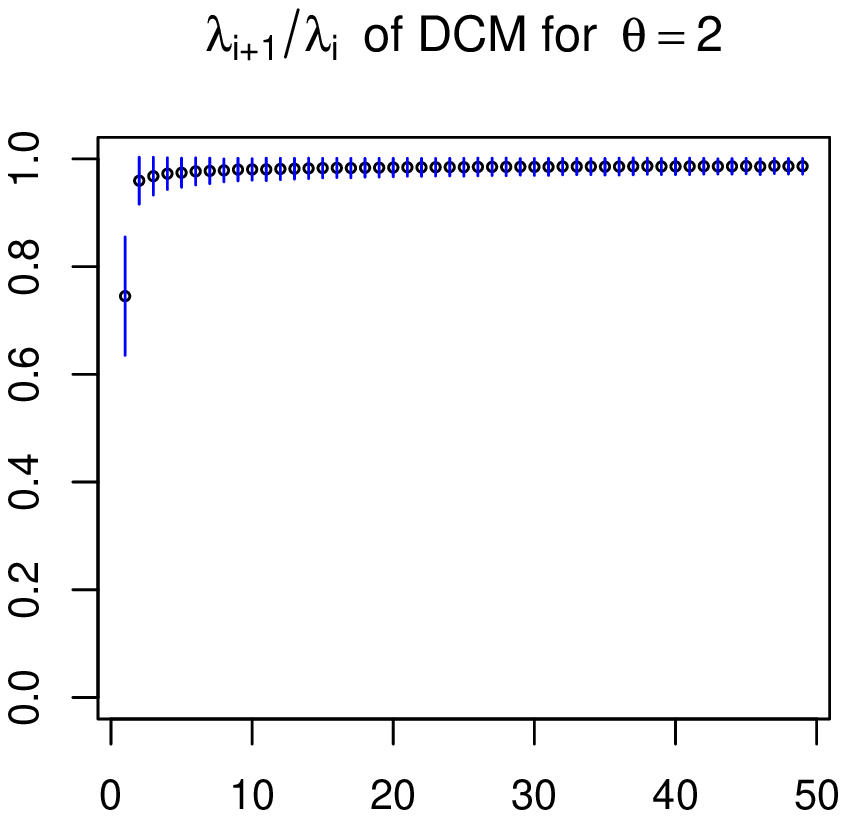}
	\end{minipage}%
	\begin{minipage}[t]{0.33\linewidth}
		\includegraphics[width=2in]{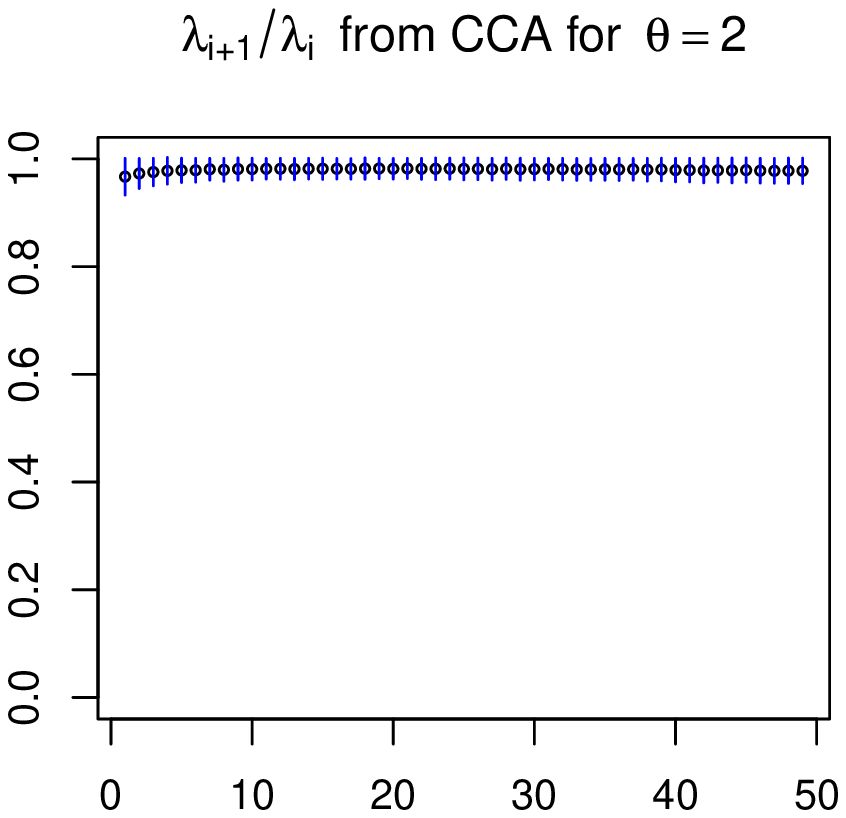}
	\end{minipage}
	\begin{minipage}[t]{0.33\linewidth}
		\includegraphics[width=2in]{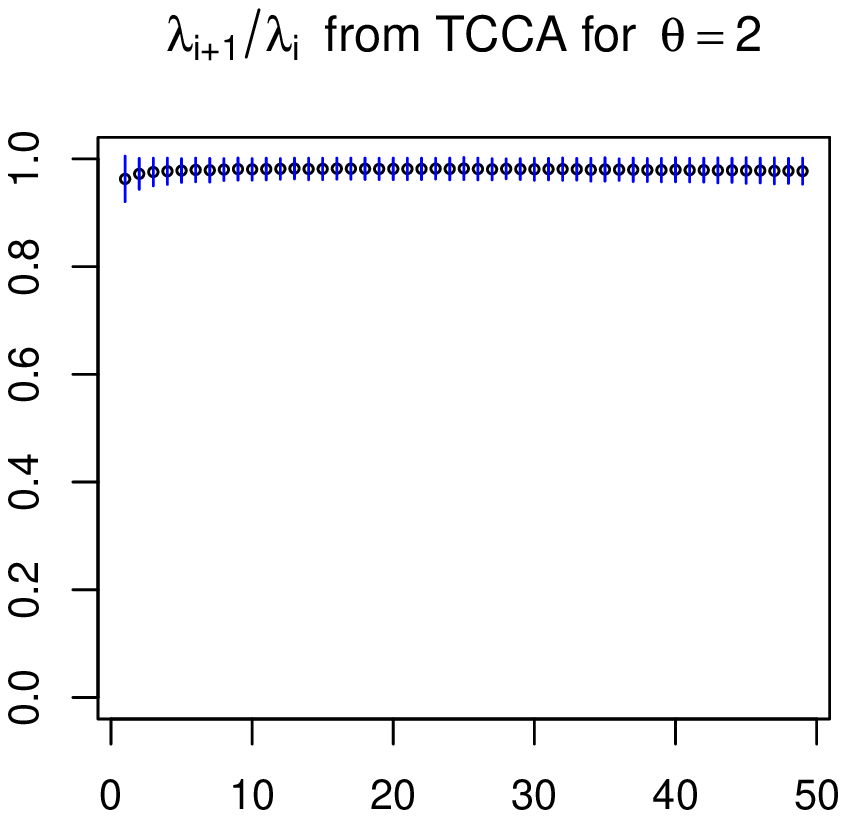}
	\end{minipage}
	\begin{minipage}[t]{0.33\linewidth}
		\includegraphics[width=2in]{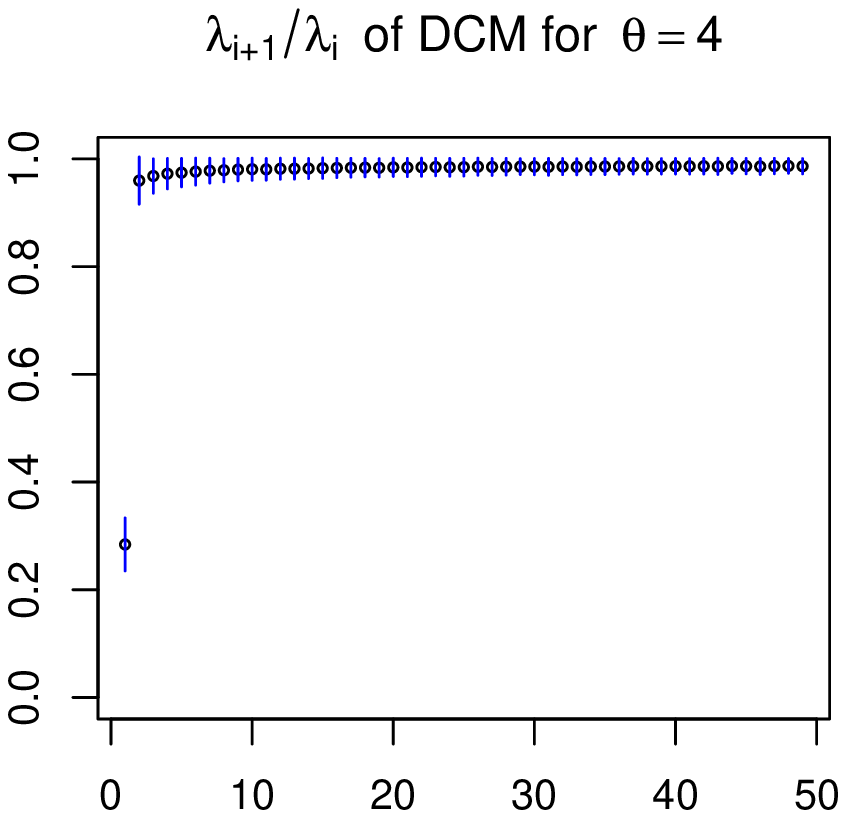}
	\end{minipage}%
	\begin{minipage}[t]{0.33\linewidth}
		\includegraphics[width=2in]{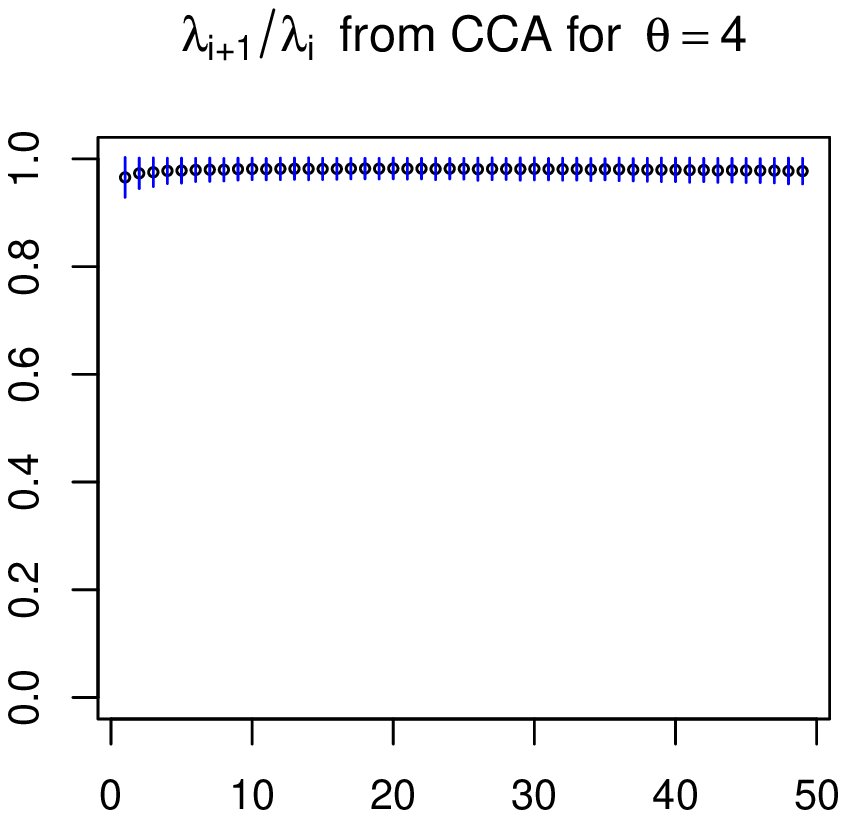}
	\end{minipage}
	\begin{minipage}[t]{0.33\linewidth}
		\includegraphics[width=2in]{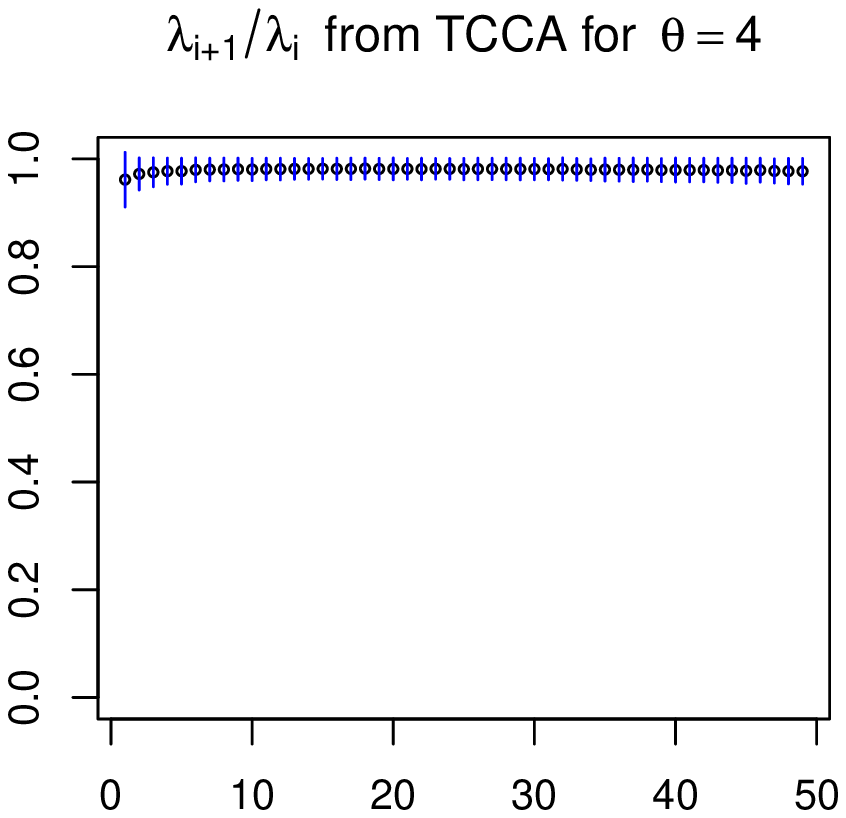}
	\end{minipage}
		\begin{minipage}[t]{0.33\linewidth}
		\includegraphics[width=2in]{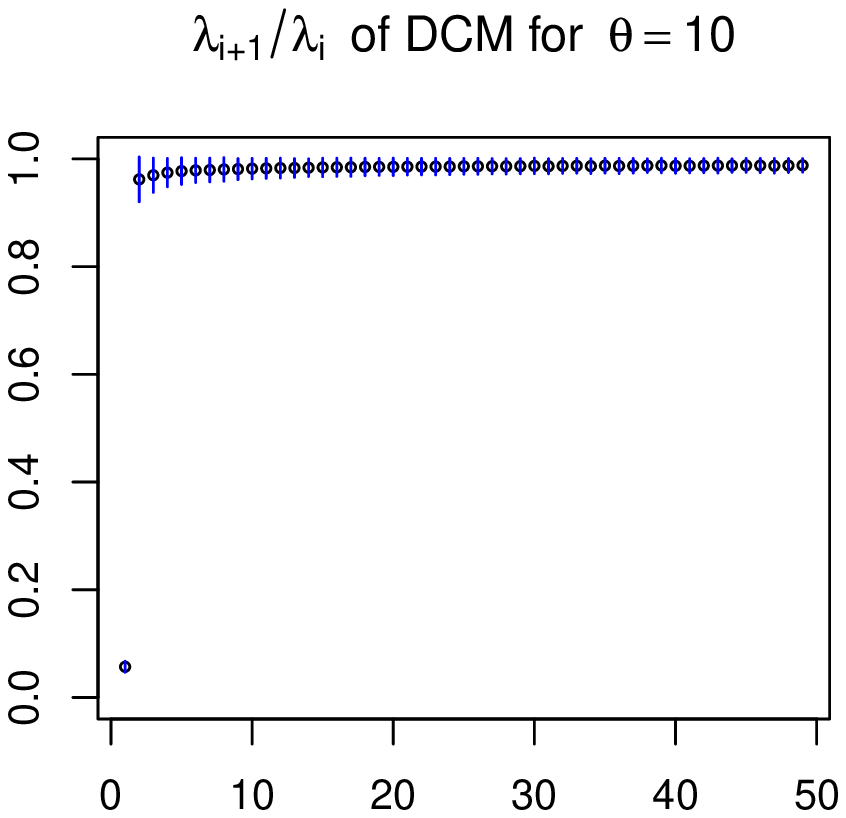}
	\end{minipage}%
	\begin{minipage}[t]{0.33\linewidth}
		\includegraphics[width=2in]{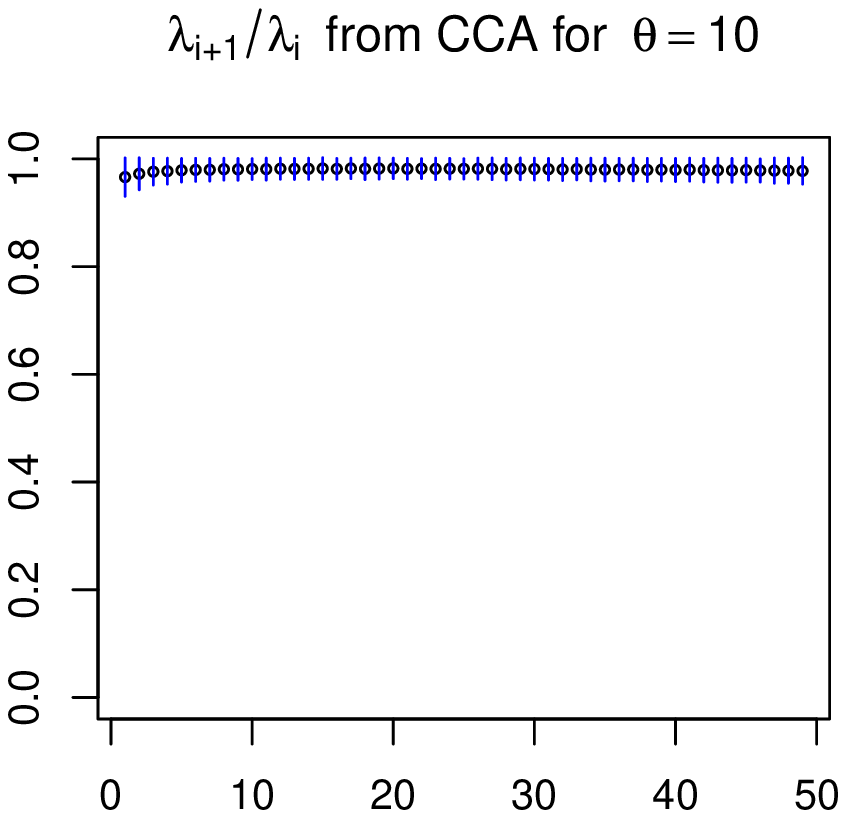}
	\end{minipage}
	\begin{minipage}[t]{0.33\linewidth}
		\includegraphics[width=2in]{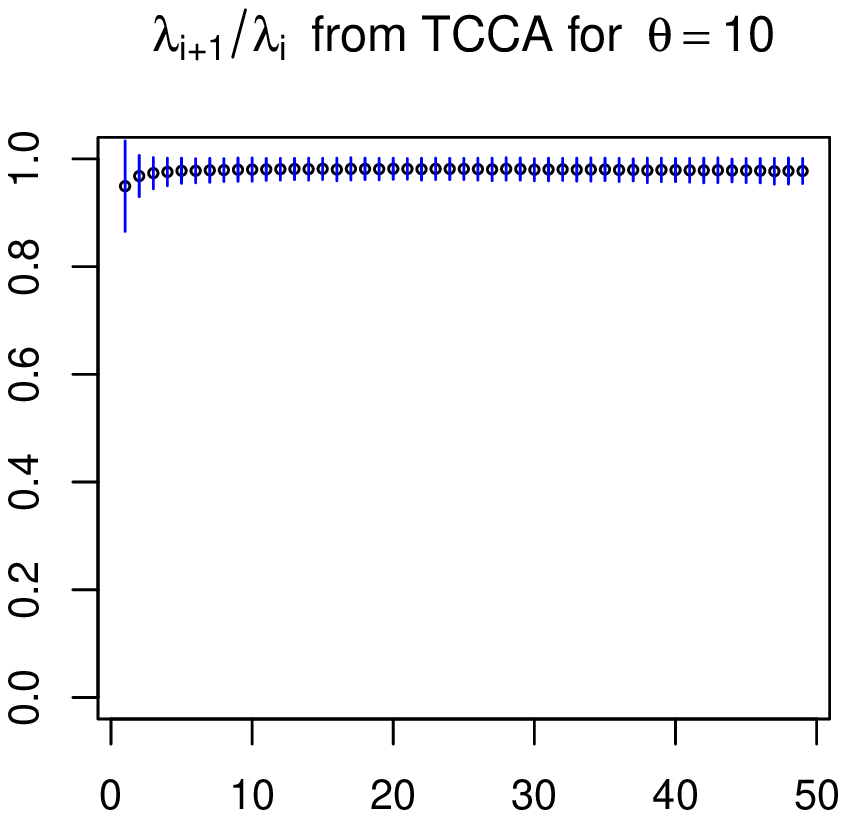}
	\end{minipage}
	\setlength{\abovecaptionskip}{-1.2cm} 
	\caption {The average of the sequences  of sample ratios $\{\lambda_{n, i+1}/ \lambda_{n, i}\}$ from 1000 independent replications under Model 5 with $\theta=2,4$. Plots on the left panel are based on the DCM $\bS_{xz}$, while those on the middle and right panels are based on CCA and TCCA, respectively. The dimensions are $(p,q,n)=(100,200,1000)$.
	}\label{rs}
\end{figure}

{
Lastly, we compare the performance of using the largest eigenvalues of our DCM model with
high-dimensional CCA \citep{yang15, B18} for detecting dependence between two
groups of random samples.  As is well known,  a direct application of CCA often fails for the detection when the two sample sets are dependent but uncorrelated. It is thus suggested  in \cite{yang15} to transform the data in a suitable way before applying CCA if one has some prior knowledge of the dependence structure. We refer to this variant of CCA  as TCCA in the following.

Model 5 is employed in this experiment. 
The parameter settings are $\theta=2,4, 10$ and
$(p,q,n)=(100,200,1000)$. For the TCCA method, we use
the exponential function $f(x)=e^x$ to transform each coordinate of the sample vectors and then conduct the CCA procedure. In this way, the two sets of transformed data are indeed linearly correlated.

Histograms of the bulk eigenvalues and the
largest eigenvalue are plotted in Figure \ref{spike-cca}. On the left
panel, the eigenvalues are from the DCM $\bS_{xz}$. One may see that the empirical SD of the bulk eigenvalues (black strips) is perfectly predicted by its LSD density curve (red lines). Moreover, the largest eigenvalues (blue strips) are centered at $\lambda=69.83, 187.5$, and $1041.5$ (blue lines) for $\theta=2, 4$ and $10$, respectively, which are clearly separated from the bulks. Similar statistics from CCA are shown on the middle panel. It demonstrates that the largest eigenvalues for $\theta= 2, 4, 10$ are all centered at $\lambda=0.49$ which is smaller than the right edge point $\lambda_+=0.5$ of the LSD. Results from TCCA are plotted on the right panel, where the largest eigenvalues are centered at $\lambda=0.49, 0.50, 0.52$ for $\theta= 2, 4, 10$, respectively.
On the other hand, Figure \ref{rs} reports the sequences of sample ratios $\{\lambda_{n, i+1}/ \lambda_{n, i}\}$ with $\pm 2$ standard deviations. For $\theta=2, 4, 10$, the first ratio $\{\lambda_{n, 2}/ \lambda_{n, 1}\}$ from the DCM model is well separated from the rest ones while those from  the CCA and TCCA models have no clear separation.
Therefore,
the
non-linear correlation between $\xx$ and $\z$ can be entirely captured
by the DCM model while CCA and TCCA will both fail to distinguish it efficiently.
We note that for the TCCA method, it indeed has some potential for the detection as one may observe that, on average, the largest eigenvalue from TCCA can surpass the right edge limit $0.5$ of the LSD as the parameter $\theta$ increases.  However its power is weak compared with our proposed method for the studied cases. Some other transforms are also tested under the same settings, such as polynomial functions, Box-Cox transforms, and  trigonometric functions. Their performance is either comparable with or less superior to the exponential function.

}

\subsection{A consistent estimator for the order of finite-rank dependence}
\label{estimating-m}

Assume that among the $m$ dependence strengths $(\theta_k)_{1\le k\le
  m}:=\btheta$, there are $m_0$ strengths above the critical
value $\theta_0$ given in \eqref{theta0}.
According to Corollary~\ref{co-spike}, the $m_0$ largest eigenvalues
$\lambda_{n, k}$ of the DCM $\bS_{xz}$
will converge almost surely to $m_0$  limits $\lambda_k$, $1\le k\le m_0$,
which are outside the
support of $F$ and given in \eqref{outlier}.
Meanwhile, the  following  eigenvalues of any given number, say $s$,
$\lambda_{n, m_0+1}, \ldots, \lambda_{n, m_0+s}$
will all converge to the right edge $\lambda_+$ of the LSD $F$.
 The rank $m_0$ corresponds to the detectable rank of the weak
  dependence considered here. In a sense, the remaining $m-m_0$
  dependence strengths $\{\theta_{m_0+1},\ldots,\theta_m\}$ below the critical
  value $\theta_0$ are too weak for detection.
Following a popular ratio estimator for the number of factors or
spikes developed in \cite{onatski10} and \cite{li17}, we introduce   a
consistent estimator for the detectable  dependence  rank $m_0$ in the
model~\eqref{zx} as follows.  Note that we have for $j=1,\ldots,m_0$,
the ratios
$
  {\lambda_{n,j+1}}/{\lambda_{n,j}}
$
will converge almost surely to a number in $(0,1)$, while for $j\ge m_0+1$
these ratios will converge to 1.
Let  $0< d_n<1$  be a  sequence of positive and vanishing constants
and consider  the following
estimator for the  dependence rank  $m_0$:
\begin{equation*}
  \hat{m}_0=
  \{ \text{first}~ j\ge 1~ ~\text{such that }~
  \frac{\lambda_{n,j+1}}{\lambda_{n,j}}   > 1-d_n\}-1.
\end{equation*}
Under the conditions similar to Theorem 3.1 in \cite{li17}, one can
show that $\hat{m}_0$ will converge to $m_0$ almost surely.

It remains
to set up an appropriate value  for the tuning parameter $d_n$.
Theoretically
any vanishing sequence $d_n\to 0$ is sufficient for the
consistence of $\hat m_0$. Here we follow the calibration proposed in
\cite{li17}.  Precisely,
we find empirically $q_{n,p,q,0.5\%}$, the lower 0.5\% quantile of
$\displaystyle  n^{2/3}\big(\nu_{2}/\nu_{1}-1\big)$
where $\nu_1$ and $\nu_2$ are the top two sample eigenvalues of the
distance covariance matrix $\bS_{xy}$ under the null model with
$\xx\sim  N({\bf 0},{\bf I}_p)$ and
$\y\sim  N({\bf 0},{\bf I}_q)$.
Then we set
$d_n = n^{-2/3} |q_{n,p,q,0.5\%}| $.
Note that $d_n$ vanishes at rate
$n^{-2/3}$. This tuned value of $d_n$ is used for all
the simulation experiments in this section.

We now examine the performance of $\hat m_0$ in finite sample
situations. Models 5 and 6 are adopted again when generating samples
of $\xx$ and $\y$. Under Model 5, we take $m=3$ and $\btheta=(4,3,2)$.
The critical value  $\theta_0$ is $1.2$, and thus the detectable  dependence rank is $m_0=3$. Under Model 6, we take $m=4$ and $\btheta=(4,3,2,1)$. In this case $\theta_0=2.5$ and $m_0=2$. Frequencies of $\hat m_0$ are calculated from 1000 independent replications under the two models with the sample size $n$ ranging from 100 to 1600. The results are shown in Tables \ref{mhat5} and  \ref{mhat6}, which verify the convergence of the proposed estimator.

\begin{table}[H]
\centering
  \caption{Frequencies  of $\hat m_0$ under Model 5
    with $\btheta=(4,3,2)$ and $m_0=3$ from 1000 independent replications. The dimensional settings are
    $c_{n1}=0.1, c_{n2}=0.2$ and $n$ ranging from 100 to 1600. \label{mhat5}}
  \begin{tabular}{cccccccc}
      \hline
		&$\hat{m}_0=0$&$\hat{m}_0=1$&$\hat{m}_0=2$&$\hat{m}_0=3$&$\hat{m}_0=4$\\
    \hline
    $n=100$ &0.045 &0.649 &0.293 &0.013 &0\\
$n=200$ &0     &0.144 &0.676 &0.176 &0.004 \\
$n=400$ &0     &0.020 &0.406 &0.561 &0.013  \\
$n=800$ &0     &0     &0.057 &0.942 &0.001   \\
$n=1600$&0     &0     &0     &0.995 &0.005   \\
\hline
	\end{tabular}
\end{table}

\begin{table}[H]
\centering
  \caption{Frequencies  of $\hat m_0$ under Model 6
    with $\btheta=(4,3,2,1)$ and $m_0=2$ from 1000 independent replications.  The dimensional settings are
    $c_{n1}=1, c_{n2}=2$ and $n$ ranging from 100 to 1600. \label{mhat6}}
  \begin{tabular}{cccccccc}
      \hline
    &$\hat{m}_0=0$&$\hat{m}_0=1$&$\hat{m}_0=2$&$\hat{m}_0=3$\\
    \hline
    $n=100$ &0.122 &0.743 &0.135 &0\\
    $n=200$ &0.016 &0.625 &0.357 &0.002   \\
    $n=400$ &0     &0.409 &0.584 &0.007   \\
    $n=800$ &0     &0.179 &0.813 &0.008   \\
    $n=1600$&0     &0.039 &0.953 &0.008   \\
    \hline
  \end{tabular}
\end{table}

\section*{Supplementary Materials}
A  supplementary to this article   contains some additional technical tools
used in this paper  and  the proofs of Theorem \ref{connection}, Theorem \ref{th:Tn},  Theorem
\ref{lsd}, Theorem \ref{lsd-spike} and Theorem  \ref{th-spike}.

\par
\section*{Acknowledgements}

{  The authors are grateful to Prof. Xiaofeng Shao for important
 discussions from which this research originated. }
 Weiming Li's research is partially supported by the NSFC (No. 11971293) and the Program of IRTSHUFE.    Qinwen Wang acknowledges support from a NSFC
 Grant (No.\ 11801085) and the Shanghai Sailing Program (No.\
 18YF1401500).
 Jianfeng Yao's  research is partially supported by a Hong Kong SAR RGC
 Grant (GRF 17308920).

\par


\bibhang=1.7pc
\bibsep=2pt
\fontsize{9}{14pt plus.8pt minus .6pt}\selectfont
\renewcommand\bibname{\large \bf References}
\expandafter\ifx\csname
natexlab\endcsname\relax\def\natexlab#1{#1}\fi
\expandafter\ifx\csname url\endcsname\relax
  \def\url#1{\texttt{#1}}\fi
\expandafter\ifx\csname urlprefix\endcsname\relax\def\urlprefix{URL}\fi


\vskip .65cm
\noindent
School of Statistics and Management, Shanghai University of Finance and    Economics
\vskip 2pt
\noindent
E-mail: \texttt{\small li.weiming@shufe.edu.cn}
\vskip 2pt

\noindent
School of Data Science,     Fudan University
\vskip 2pt
\noindent
E-mail:  \texttt{\small wqw@fudan.edu.cn}
\vskip 2pt

\noindent
Department of Statistics and Actuarial Science, The University of Hong Kong
\vskip 2pt
\noindent
E-mail:  \texttt{\small jeffyao@hku.hk}
\vskip 2pt

\newpage

\appendix
\newcommand{\G}{{\mathbf G}}
\newcommand{\bH}{{\mathbf H}}

\renewcommand{\baselinestretch}{2}

\markright{ \hbox{\footnotesize\rm Supplementary material
}\hfill\\[-13pt]
\hbox{\footnotesize\rm
}\hfill }

\markboth{\hfill{\footnotesize\rm  W. Li, Q. Wang and J. Yao.} \hfill}
{\hfill {\footnotesize\rm High dimensional distance covariance matrix} \hfill}

\renewcommand{\thefootnote}{}
$\ $\par \fontsize{12}{14pt plus.8pt minus .6pt}\selectfont


 \centerline{\large\bf Eigenvalue distribution of a high-dimensional}
\vspace{2pt}
 \centerline{\large\bf distance  covariance matrix with application}
\vspace{.25cm}
 \centerline{Weiming Li, Qinwen Wang and Jianfeng Yao}
 \vspace{.4cm}
\centerline{\it Shanghai University of Finance and    Economics,}
\centerline{\it Fudan University and The University of Hong Kong}
 
\vspace{.55cm}
 \centerline{\bf Supplementary Material}
\vspace{.55cm}
\fontsize{9}{11.5pt plus.8pt minus .6pt}\selectfont
\noindent

This supplementary material contains some additional technical tools
  and  the proofs of Theorem \ref{connection}, Theorem \ref{th:Tn},  Theorem
\ref{lsd}, Theorem \ref{lsd-spike} and Theorem  \ref{th-spike} of the
main paper.   \red{Throughout this  supplementary material, $||\cdot||$ denotes the Euclidean norm for vectors, the spectral norm for matrices and the supremum norm for functions, respectively. ${\mathbb C}^+$ and ${\mathbb C}^-$ are referred as  the upper and lower  half
  complex plane  (real axis excluded). $K$ is used to denote some constant that can vary from place to place.

 }

\par

\fontsize{12}{14pt plus.8pt minus .6pt}\selectfont

\section{Technical tools}

\begin{lemma}\label{kernel}[\citet{el10}]
	Consider the $n \times n$ kernel random matrix $\M$ with entries
	$$\M_{i, j}=f\left(\frac{\|\xx_i-\xx_j\|^2_2}{p}\right).$$
	Let us call $\Ps$ the vector with $i$-th entry $\Ps_i=\|\xx_i\|^2_2/p-\tau/2$, where $\tau=2\tr(\Sig_p)/p$. We assume that:\\
	(a) $n\asymp p$, that is, $n/p$ and $p/n$ remain bounded as $p \to \infty$.\\
	(b) $\Sig_p$ is a positive semi-definite $p\times p$ matrix, and $\|\Sig_p\|=\sigma_1(\Sig_p)$ remains bounded in $p$, that is, there exists $K>0$, such that $\sigma_1(\Sig_p)\leq K$, for all $p$.\\
	(c) There exists $\ell \in \mathbb R$ such that $\lim_{p \to \infty}\tr(\Sig_p)/p=\ell$.\\
	(d) $\X=(\xx_1,\ldots,\xx_n)$ and $\xx_i=\Sig^{1/2}_p\w_i$ for $i=1,\ldots,n$.\\
	(e) The entries of $\w_i$, a $p$-dimensional random vector, are i.i.d. Also, denoting by $w_{ik}$ the $k$th entry of $\w_i$, we assume that $\E(w_{ik}) = 0$, $\var(w_{ik}) = 1$ and $\E(|w_{ik}|^{5+\varepsilon}) <\infty$ for some $\varepsilon> 0$. \\
	(f) $f$ is $C^3$ in a neighborhood of $\tau$.\\
	Then $\M$ can be approximated consistently in operator norm (and in probability) by the matrix $\widetilde{\M}$, defined by
	\begin{align*}
	\widetilde{\M}&=f(\tau)\1 \1'+{f}'(\tau)\left[\1 \Ps'+\Ps {\1'}-2\frac{\X'\X}{p}\right]\\
	&\quad +\frac{f''(\tau)}{2}\left[\1(\Ps \circ \Ps)'+(\Ps \circ \Ps){\1'}+2\Ps \Ps'+4\frac{\tr(\Sig^2_p)}{p^2}\1\1'\right]+v_p {\bf I}_n,\\
	v_p&=f(0)+\tau f'(\tau)-f(\tau).
	\end{align*}
	In other words,
	$$||\M-\widetilde{\M}|| \to 0,\quad\text{in~ probability}.$$	
\end{lemma}

\begin{lemma}
	\label{rank-ineq}[\citet{BSbook}]
	Let $\A$ and $\B$ be two $n\times n$ Hermitian matrices. Then,
	$$
	||F^{\A}-F^{\B}||\leq \frac{1}{n}{\rm rank}(\A-\B)\quad \text{and}\quad L^3(F^\A,F^\B)\leq \frac{1}{n}\tr[(\A-\B)(\A-\B)^*],
	$$
	where $L(F,G)$ stands for the L\'evy distance between the distribution functions $F$ and $G$.
\end{lemma}

\begin{lemma}\label{cth}
	Let $f$: $\mathbb R^n \to \mathbb R$ be any function of thrice differentiable in each argument.
	Let also $\xx=(x_1,\ldots,x_n)'$ and $\y=(y_1,\ldots,y_n)'$ be two random vectors in $\mathbb R^n$ with i.i.d.\ elements, respectively, and set $U=f(\xx)$ and $V=f(\y)$. If $$\gamma=\max \{\E|x_i|^3, \E|y_i|^3, 1\leq i\leq n\}<\infty,$$
	then for any thrice differentiable $g: \mathbb R \to \mathbb R$ and any $K>0$,
	\begin{align*}
	|\E g(U)-\E g(V)|\leq 2C_2(g)\gamma n\lambda_3(f),
	\end{align*}
	where $\lambda_3(f)=\sup\left \{|\partial ^k_i f(\z)|^{3/k}:  \z=(z_\ell), z_\ell\in \{x_\ell,y_\ell\}, 1\leq i\leq n, 1\leq k\leq 3 \right\}$ and $C_2(g)=\frac 16\|g'\|_{\infty}+\frac 12\|g''\|_{\infty}+\frac 16 \|g'''\|_{\infty}$.
\end{lemma}
This lemma follows directly from Corollary 1.2 in \cite{ch08} and its proof.

\section{Proofs}\label{proof}
At the beginning of this section, we first recall some notations for easy reading.
\begin{align*}
&\V_x=\left(\frac{\|{\bf x}_k-\bx_\ell\|}{\sqrt{p}}\right),\quad
\V_y=\left(\frac{\|{\bf y}_k-\by_\ell\|}{\sqrt q}\right),\quad {\bf P}_n=\I_n-\frac 1n {\bf 1}_n{\bf 1}_n',\\
  &\gamma_x= \frac 1p \tr {\Sig_x},~\quad \gamma_y= \frac 1q \tr \Sig_y,\quad \kappa_{x}=\frac{1}{pn}\sum_{i=1}^n||\xx_i||^2, \quad 
  \kappa_{y}=\frac{1}{qn}\sum_{i=1}^n||\y_i||^2, \nonumber\\
     &\A_n= \frac{1}{p}\X'\X+\gamma_x{\bf I}_n,~\quad  \C_n=\frac{1}{q}\Y'\Y+\gamma_y{\bf I}_n,\quad \B_n=\A_n^{\frac12}\C_n\A_n^{\frac12},\\
      & {\D_x}=\frac{1}{p}{\X'\X}+\kappa_{x} {\I_n},\quad {\D_y}= \frac{1}{q} {\Y'\Y+\kappa_{y} \I_n}, \quad\D_z=\frac{1}{q}{\Z'\Z}+\kappa_{z}{\bf I}_n,\\
      &  \bS_{xy}={\bf P}_n\D_x{\bf P}_n\D_y{\bf P}_n,\quad \bS_{xz}=\bP_n \D_x \bP_n\D_z\bP_n.
\end{align*}

\subsection{Proof of Theorem \ref{connection}} 

The squared sample distance covariance
$\mV^2_n(\bx,\by)$ in \eqref{eq:Vn} can be expressed as an inner product between the two matrices $\bP_n\V_x\bP_n$ and $\bP_n\V_y\bP_n$, that is,
$$ 
\mV_n^2(\bx, \by) = \frac{\sqrt{pq}}{n^2} \tr\bP_n\V_x \bP_n \V_y \bP_n.
$$
Notice that the matrices $\V_x$ and $\V_y$ are exactly the
Euclidean distance kernel matrices discussed in \citet{el10}
with  kernel   function   $f(x)=\sqrt x$.
Applying their main theorem (see Lemma~\ref{kernel}), the matrix
\begin{align}\label{vnde}
 \bP_n\V_x \bP_n \V_y \bP_n
\end{align}
can be approximated by a simplified random matrix $\V_n$ such that as $(n,p,q)$ tend to infinity,
 \begin{align}\label{vnvn1}
 \left \|\V_n -  \bP_n\V_x \bP_n \V_y \bP_n\right\|\to 0
  \end{align}
   in probability, where
\begin{align}\label{nnab}
  \V_{n}\triangleq \frac{1}{2\sqrt{\gamma_x \gamma_y}}\bP_n \left(\A_n+\frac{1}{8\gamma_x}\Ps_x\Ps'_x\right)\bP_n \left(\C_n+\frac{1}{8\gamma_y}\Ps_y\Ps'_y\right)\bP_n,
\end{align}
in which
\begin{align}
  &{\boldsymbol\psi_x}=\frac1p
    \left(\begin{array}{c}
            \|\bx_1\|^2-\tr\Sig_x\\
            \vdots\\
            \|\bx_n\|^2-\tr\Sig_x
          \end{array}\right)~~\text{and}~~
  \Ps_y=\frac1q
  \left(\begin{array}{c}
          \|\by_1\|^2-\tr\Sig_y\\
          \vdots\\
          \|\by_n\|^2-\tr\Sig_y
        \end{array}\right).\nonumber
\end{align}

Then we replace the two traces $\gamma_x$ and $\gamma_y$ in $\A_{n}$ and $\C_n$ with their unbiased sample counterparts $\kappa_{x}$ and $\kappa_{y}$, respectively,
which does not affect the convergence in \eqref{vnvn1}.
Finally in \eqref{nnab}, 
by  removing the two  rank-one matrices
$(8\gamma_x)^{-1}\Ps_x\Ps^T_x$ and $(8\gamma_y)^{-1}\Ps_y\Ps^T_y$ (which have bounded spectral norm, almost surely), we get the conclusion of the theorem.
The proof is thus complete.

\subsection{Proof of Theorem \ref{th:Tn}} 
{ Recall the approximation from Theorem \ref{connection},
\begin{align*}
  \mathcal{V}^2_n( \bx,\by)=\frac{1}{2n^2}\sqrt{\frac{pq}{\gamma_x\gamma_y}}\tr \bS_{xy}+o_p(1)
\end{align*}
and notice that
\begin{align*}
  \frac{1}{n}\tr(\bS_{xy})&=\frac{1}{npq}\tr(\bP_n\X'\X\bP_n\Y'\Y\bP_n)+\frac{\kappa_y}{np}\tr(\bP_n\X'\X)+\frac{\kappa_x}{nq}\tr(\bP_n\Y'\Y)\\
  &\quad+\frac{n-1}{n}\kappa_x\kappa_y\\
		&=\frac{1}{npq}\tr(\X'\X\Y'\Y)+3\gamma_x\gamma_y+o_{a.s}(1).
		\end{align*}
		Moreover, from Equation (21) in \citet{LY18} and the independence between $\X$ and $\Y$,
		\begin{align*}
		\frac{1}{npq}\tr(\X'\X\Y'\Y)&=\frac{1}{p}\tr(\Sig_x)\frac{1}{q}\tr(\Sig_y)+o_{a.s}(1).
		\end{align*}
		Collecting the above results yields
		\begin{align*}
		\mathcal{V}^2_n( \bx,\by)=2\sqrt{c_{n1}c_{n2}\gamma_x\gamma_y}+o_p(1).
		\end{align*}
	    On the other hand, applying Lemma \ref{kernel}, we have
		\begin{align*}
		\frac{1}{n}S_{2,n}=&\frac{1}{2n}\sqrt{\frac{pq}{\gamma_x\gamma_y}}\left(\frac{1}{n^2}\1'\D_x\1-2\gamma_x\right)\left(\frac{1}{n^2}\1'\D_y\1-2\gamma_y\right)+o_p\left(1\right)\nonumber\\
		=&2\sqrt{c_{n1}c_{n2}\gamma_x\gamma_y}+o_p(1).
		\end{align*}
		Therefore, the statistic $T_n=n\mathcal{V}^2_n( \bx,\by)/S_{2,n}$ converges to 1 in probability.
	The proof is complete.

\subsection{Proof of Theorem \ref{lsd}}

The strategy of the proof is as follows. First, we prove the theorem under Gaussian assumption. By virtue of rotation invariance property of Gaussian vectors, we may treat the two population covariance matrices $\Sig_x$ and $\Sig_y$ as diagonal ones, which can simplify the proof dramatically. Second,  applying Lindeberg's replacement trick provided in \cite{ch08}, we will remove the Gaussian assumption and show that the theorem still holds true for general distributions if the atoms $(w_{ij})$ have finite fourth moment, as stated in our  Assumption (b).
\medskip

\noindent {\em Gaussian case:}
First, we have
\begin{align}\label{k-r}
  |\kappa_x-\gamma_x|\xrightarrow{a.s.}0\quad\text{and}\quad|\kappa_y-\gamma_y|\xrightarrow{a.s.}0,
\end{align}
as $(n,p,q)$ tend to $\infty$.
>From Lemma \ref{rank-ineq} and \eqref{k-r}, we get
\begin{align*}
  L^3(F^{\bS_{xy}},F^{\B_n})\xrightarrow{a.s.}0.
\end{align*}		
Hence, the matrices $\bS_{xy}$ and $\B_n$ share the same limiting spectral distribution
and thus we only focus on the convergence of $F^{\B_n}$. We first derive its limit conditioning on the sequence $(\A_n)$.
Then the result holds unconditionally if the limit is independent of  $(\A_n)$. Following standard strategies from random matrix theory, letting $s_{\B_n}(z)$ be the Stieltjes transform of $F^{\B_n}$, the convergence of $F^{\B_n}$ can be established through three steps:
\begin{itemize}
\item[] {\em Step 1}: For any fixed $z\in \mathbb C^+$, $s_{\B_n}(z)-\E s_{\B_n}(z)\to0$, almost surely.
\item[] {\em Step 2}: For any fixed $z\in \mathbb C^+$, $\E s_{\B_n}(z)\to s(z)$ with $s(z)$ satisfies the equations in \eqref{lsd}.
\item[] {\em Step 3}: The uniqueness of the solution $s(z)$ to \eqref{lsd-sys} on the set 	\eqref{set}.
\end{itemize}
\smallskip

\noindent{\em Step~1. Almost sure convergence of $s_{\B_n}(z)-\E s_{\B_n}(z)$.}
\medskip

We assume $\Sig_y$ is diagonal, having the form
$$\Sig_y={\rm Diag}(\tau_1,\ldots, \tau_q).$$
By this and notations
\begin{align*}
  \br_k=\frac{1}{\sqrt{q}}\A_n^{1/2}(w_{p+k,1},\ldots,w_{p+k,n})',\quad k=1,\ldots,q,
\end{align*}
the matrix $\B_n$ can be expressed as
\begin{align}\label{bn}
  \B_n=\gamma_y\A_n+\sum_{k=1}^q\tau_k\br_k\br_k'.
\end{align}
It's ``leave-one-out" version is denoted by $\B_{k,n}=\B_n-\tau_k\br_k\br_k'$, $k=1,\ldots,q$.
Let $\E_0(\cdot)$ be expectation and $\E_k(\cdot)$ be conditional expectation given $\br_{1},\ldots, \br_k$. From the martingale decomposition and the identity
\begin{align}\label{equ0}
  \br_k'(\B_n-z\I_n)^{-1}=\frac{\br_k'(\B_{k,n}-z\I_n)^{-1}}{1+\tau_k\br_k'(\B_{k,n}-z\I_n)^{-1}\br_k},	
\end{align}
we have
\begin{align}
  s_{\B_n}(z)-\E s_{\B_n}(z)=&\frac{1}{n}\sum_{k=1}^q(\E_k-\E_{k-1})\left[\tr(\B_n-z\I_n)^{-1}-\tr(\B_{k,n}-z\I_n)^{-1}\right]\nonumber\\
  =&-\frac{1}{n}\sum_{k=1}^q(\E_k-\E_{k-1})\frac{\tau_k\br_k'(\B_{k,n}-z\I_n)^{-2}\br_k}{1+\tau_k\br_k'(\B_{k,n}-z\I_n)^{-1}\br_k}.\label{mds}
\end{align}
Similar to the arguments on pages 435-436 of \cite{BZ08}, the summands in \eqref{mds} form a bounded martingale difference sequence, and hence $s_{\B_n}(z)-\E s_{\B_n}(z)\to0$, almost surely.

\medskip

\noindent{\em Step~2. Convergence of $\E s_{\B_n}(z)$.}
\medskip

Let $s_{\A_n}(z)$ be the Stieltjes transform of $F^{\A_n}$. From \cite{S95},  $s_{\A_n}(z)$ converges almost surely to $s_\A(z)$, which satisfies
\begin{align}\label{sa}
  z=-\frac{1}{s_\A(z)}+\int t+\frac{t}{1+tc_1^{-1}s_\A(z)}dH_x(t).
\end{align}
Define two functions $w_n(z)$ and $m_n(z)$ as
\begin{align}\label{wmn}
  w_n(z)=\frac{1}{n}\E\tr(\B_n-zI_n)^{-1}\A_n \quad\text{and}\quad m_n(z)=\gamma_y+\frac{1}{q}\sum_{k=1}^q\frac{\tau_k}{1+\tau_kc_{n2}^{-1} w_n(z)}.
\end{align}
We first show that
\begin{align}\label{lsd:c1}
  m_n^{-1}(z)s_{\A_n}\left[zm_n^{-1}(z)\right]-\E s_{\B_n}(z)\to 0,\quad n\to\infty.
\end{align}
In fact, applying the identity \eqref{equ0}, we have
\begin{align*}
  &\frac{1}{n}\tr\left[m_n(z)\A_n-z\I_n\right]^{-1}-\frac{1}{n}\tr(\B_n-z\I_n)^{-1}\\
  =&\frac{1}{n}\tr\left[m_n(z)\A_n-z\I_n\right]^{-1}\left(\sum_{k=1}^q\tau_k\br_k\br_k'-(m_n(z)-\gamma_y)\A_n\right)(\B_n-z\I_n)^{-1}\no
     =&\frac{1}{n}\sum_{k=1}^q\frac{\tau_k\br_k'(\B_{k,n}-zI_n)^{-1}\left[(m_n(z)\A_n-z\I_n\right]^{-1}\br_k}{1+\tau_k\br_k'(\B_{k,n}-z\I_n)^{-1}\br_k}\no
      &-\frac{m_n(z)-\gamma_y}{n}\tr \left[m_n(z)\A_n-z\I_n\right]^{-1}\A_n(\B_n-z\I_n)^{-1}\no
        =&\frac{1}{n}\sum_{k=1}^q\frac{\tau_kd_k}{1+\tau_kc_{n2}^{-1} w_n(z)},\nonumber
\end{align*}
where
\begin{align*}
  d_k=&\frac{1+\tau_kc_{n2}^{-1}w_{n}(z)}{1+\tau_k\br_k'(\B_{k,n}-z\I_n)^{-1}\br_k}\br_k'(\B_{k,n}-z\I_n)^{-1}\left[m_n(z)\A_n-z\I_n\right]^{-1}\br_k\\
      &-\frac{1}{q}\tr \left[m_n(z)\A_n-z\I_n\right]^{-1}\A_n(\B_n-z\I_n)^{-1}.
\end{align*}
Following similar arguments on pages 85-87 of \cite{BSbook}, one may obtain
$$
\max_k \E (d_k)\to0.
$$
This result together with the fact
$$
\inf_n|1+\tau_kc_{n2}^{-1}w_n(z)|\geq \inf_n\tau_kc_{n2}^{-1}|\Im(w_n(z))|>0
$$
imply the convergence in \eqref{lsd:c1}.
\medskip

We next find another link between $\E s_{\B_n}(z)$ and $w_n(z)$ by proving
\begin{align}\label{lsd:c2}
  1+z\E s_{\B_n}(z)-\gamma_yw_n(z)-\frac{1}{n}\sum_{k=1}^q\frac{\tau_kw_{n}(z)}{c_{n2}+\tau_kw_{n}(z)}\to0.
\end{align}
>From the expression of $\B_n$ in \eqref{bn} and the identity in \eqref{equ0}, we have
\begin{align}\label{e2:1}
  \I_n+z(\B_n-z\I)^{-1}=&\B_n(\B_n-z\I_n)^{-1}\no
                       =&\gamma_y\A_n(\B_n-z\I)^{-1}+\sum_{k=1}^q\tau_k\br_k\br_k'(\B_{n}-z\I_n)^{-1}\no
                          =&\gamma_y\A_n(\B_n-z\I)^{-1}+\sum_{k=1}^q\frac{\tau_k\br_k\br_k'(\B_{k,n}-z\I_n)^{-1}}{1+\tau_k\br_k'(\B_{k,n}-z\I_n)^{-1}\br_k}.
\end{align}
Taking the trace on both sides of \eqref{e2:1} and dividing by $n$, we get
\begin{align*}
  1+z\frac{1}{n}\tr (\B_n-z\I_n)^{-1}=&\gamma_y\frac{1}{n}\tr(\B_n-z\I_n)^{-1}\A_n+\frac{1}{n}\sum_{k=1}^q\frac{\tau_k\br_k'(\B_{k,n}-z\I_n)^{-1}\br_k}{1+\tau_k\br_k'(\B_{k,n}-z\I_n)^{-1}\br_k}\no
                                      =&\gamma_y\frac{1}{n}\tr(\B_n-z\I_n)^{-1}\A_n+\frac{1}{n}\sum_{k=1}^q\frac{\tau_kc_{n2}^{-1}w_{n}(z)}{1+\tau_kc_{n2}^{-1}w_{n}(z)}+\varepsilon_n,
\end{align*}
where
$$
\varepsilon_n=\frac{1}{n}\sum_{k=1}^q \frac{\tau_k[c_{n2}^{-1}w_{n}(z)-\br_k'(\B_{k,n}-z\I_n)^{-1}\br_k]}{[1+\tau_k\br_k'(\B_{k,n}-z\I_n)^{-1}\br_k][1+\tau_kc_{n2}^{-1}w_{n}(z)]}.
$$
>From the proof of (2.3) in \cite{S95}, almost surely,
$$
\inf_n\left|[1+\tau_k\br_k'(\B_{k,n}-z\I_n)^{-1}\br_k][1+\tau_kc_{n2}^{-1}w_{n}(z)]\right|>0.
$$
Moreover, following similar arguments on page 87 of \cite{BSbook}, one may get
$$
\frac{1}{n}\sum_{k=1}^q \E^{\frac{1}{2}}|c_{n2}^{-1}w_{n}(z)-\br_k'(\B_{k,n}-zI_n)^{-1}\br_k|^2\to0.
$$
Therefore $\E(\varepsilon_n)\to0$ and hence the convergence in \eqref{lsd:c2} holds.
\medskip

By considering a subsequence $\{n_k\}$ such that $w_{n_k}(z)\to w(z)$, from \eqref{sa}, \eqref{lsd:c1} and \eqref{lsd:c2}, we have
\begin{align*}
  m_{n_k}(z)&\to \int t+\frac{t}{1+tc_{n2}^{-1} w(z)}dH_{y}(t)\triangleq m(z),\\
  s_{n_k}(z)&\to\frac{1}{m(z)}s_{\A}\left(\frac{z}{m(z)}\right),\\
  zs_{n_k}(z)&\to-1+w(z)\int t+\frac{t}{1+tc_2^{-1}w(z)}dH_y(t),
\end{align*}
as $k\to\infty$.
These results demonstrate that $s_{n_k}(z)$ has a limit, say $s(z)$, which together with $(w(z),m(z)$, $s_A(z))$ satisfy the following system of equations:
\begin{align*}
\left\{
\begin{array}{l}\displaystyle
  s(z)=\frac{1}{m(z)}s_{\A}\left(\frac{z}{m(z)}\right),\\
  \displaystyle zs(z)=-1+w(z)\int t+\frac{t}{1+tc_2^{-1}w(z)}dH_y(t),\\
 \displaystyle z=-\frac{1}{s_\A(z)}+\int t+\frac{t}{1+tc_1^{-1}s_\A(z)}dH_x(t),\\
 \displaystyle m(z)=\int t+\frac{t}{1+tc_{2}^{-1}w(z)}dH_y(t).
  \end{array}
  \right.
\end{align*}
Cancelling the function $s_\A(z)$ from the above system yields an equivalent but simpler system of equations as shown in \eqref{lsd}.
Hence, the convergence of $s_n(z)$ is established if the system has a unique solution on the set \eqref{set}.

\medskip

\noindent{\em Step~3. Uniqueness of the solution to \eqref{lsd}.}
\medskip

The system of equations in \eqref{lsd} is equivalent to
\begin{align}\label{lsdt}
  \left\{
  \begin{array}{l}
    \displaystyle 1+zs=wm,\\
    \displaystyle m=\int t+\frac{t}{1+tc_2^{-1}w}dH_y(t),\\
    \displaystyle w=s\int t+\frac{t}{1+tc_1^{-1}(1+zs)w^{-1}s}dH_x(t).
  \end{array}
  \right.
\end{align}
Bringing $s=[wm-1]/z$ into the third equation in \eqref{lsdt}, we have
\begin{align}\label{lsdone}
  w=\int \frac{t}{z}\Big(wm-1\Big)+\frac{t\big(wm-1\big)}{z+c_1^{-1}tm\big(wm-1\big)}dH_x(t).
\end{align}
Now suppose the LSD $F\neq \delta_0$ and we have two solutions $(s,w,m)$ and $(\tilde s,\tilde w,\tilde m)$ to the system on the set \eqref{set}   
for a common $z \in \mathbb{C}^{+}$. Then, from \eqref{lsdt} and \eqref{lsdone}, we can obtain
\begin{align}
  &w-\tilde{w}=(wm-\tilde{w}\tilde{m})\nonumber\\
  &\quad\quad \quad\quad\times\int \left[\frac tz+\frac{tz}{\big(z+c^{-1}_1tm(wm-1)\big)\big(z+c^{-1}_1t\tilde{m}(\tilde{w}\tilde{m}-1)\big)}\right]dH_x(t)\nonumber\\ 
  &\quad\quad \quad\quad+(\tilde{m}-m)\int\frac{t^2c^{-1}_1(wm-1)(\tilde{w}\tilde{m}-1)}{(z+c^{-1}_1tm(wm-1))(z+c^{-1}_1t\tilde{m}(\tilde{w}\tilde{m}-1))}dH_x(t),\label{ww}\\
  &\tilde{m}-m=(w-\tilde{w})\int \frac{t^2c^{-1}_2}{\big(1+tc^{-1}_2w\big)\big(1+tc^{-1}_2\tilde{w}\big)}dH_y(t),\label{aa}\\
  &wm-\tilde{w}\tilde{m}=(w-\tilde{w})\int \left(t+\frac{t}{\big(1+tc^{-1}_2w\big)\big(1+tc^{-1}_2\tilde{w}\big)}\right)dH_y(t).~\label{wawa}
\end{align}
Combining \eqref{ww}-\eqref{wawa}, if $w\neq \tilde{w}$, we have
\begin{align}\label{e1}
  B_1B_2+C_1C_2=1,
\end{align}
where
\begin{align*}
  &B_1=\int \frac tz+\frac{tz}{\big(z+c^{-1}_1tm(wm-1)\big)\big(z+c^{-1}_1t\tilde{m}(\tilde{w}\tilde{m}-1)\big)}dH_x(t),\\
  &B_2=\int t+\frac{t}{\big(1+tc^{-1}_2w\big)\big(1+tc^{-1}_2\tilde{w}\big)}dH_y(t),\\
  &C_1=\int\frac{t^2c^{-1}_1(wm-1)(\tilde{w}\tilde{m}-1)}{(z+c^{-1}_1tm(wm-1))(z+c^{-1}_1t\tilde{m}(\tilde{w}\tilde{m}-1))}dH_x(t),\\
  &C_2=\int \frac{t^2c^{-1}_2}{\big(1+tc^{-1}_2w\big)\big(1+tc^{-1}_2\tilde{w}\big)}dH_y(t).
\end{align*}
By the Cauchy-Schwarz inequality, we have
\begin{align*}
  |B_1B_2|^2\leq&
                \int \Big|\frac tz\Big|+\frac{|tz|}{|z+c^{-1}_1tm(wm-1)|^2}dH_x(t) \\
                &\times \int \Big|\frac tz\Big|+\frac{|tz|}{|z+c^{-1}_1t\widetilde{m}(\widetilde{w}\widetilde{m}-1)|^2}dH_x(t)\\
              &\times \int t+\frac{t}{|1+tc^{-1}_2w|^2}dH_y(t) \int t+\frac{t}{|1+tc^{-1}_2\widetilde{w}|^2}dH_y(t)\\
  =&\int \Big|\frac tz\Big|+\frac{|tz|}{|z+c^{-1}_1tm(wm-1)|^2}dH_x(t) \int t+\frac{t}{|1+tc^{-1}_2w|^2}dH_y(t)\\
              &\times \int \Big|\frac tz\Big|+\frac{|tz|}{|z+c^{-1}_1t\widetilde{m}(\widetilde{w}\widetilde{m}-1)|^2}dH_x(t)
                \int t+\frac{t}{|1+tc^{-1}_2\widetilde{w}|^2}dH_y(t)\\
  :=& (\widetilde{B}_1\widetilde{B}_2)^2,\\
  \\
  |C_1C_2|^2\leq& \int \frac{t^2c^{-1}_1|wm-1|^2}{|z+c^{-1}_1tm(wm-1)|^2}dH_x(t)\int \frac{t^2c^{-1}_1|\widetilde{w}\widetilde{m}-1|^2}{|z+c^{-1}_1t\widetilde{m}(\widetilde{w}\widetilde{m}-1)|^2}dH_x(t)\\
              &\times \int \frac{t^2c^{-1}_2}{|1+tc^{-1}_2w|^2}dH_y(t)\frac{t^2c^{-1}_2}{|1+tc^{-1}_2\widetilde{w}|^2}dH_y(t)\nonumber\\
  =&\int \frac{t^2c^{-1}_1|wm-1|^2}{|z+c^{-1}_1tm(wm-1)|^2}dH_x(t) \int \frac{t^2c^{-1}_2}{|1+tc^{-1}_2w|^2}dH_y(t)\\
              &\times      \int \frac{t^2c^{-1}_1|\widetilde{w}\widetilde{m}-1|^2}{|z+c^{-1}_1t\widetilde{m}(\widetilde{w}\widetilde{m}-1)|^2}dH_x(t)
                \int \frac{t^2c^{-1}_2}{|1+tc^{-1}_2\widetilde{w}|^2}dH_y(t)\nonumber\\
  :=& (\widetilde C_1\widetilde C_2)^2.~\nonumber
\end{align*}
Then \eqref{e1} implies
\begin{align}\label{ine1}
  1&=|B_1B_{2}+C_1C_2|\nonumber\\
    &\leq \sqrt{(\widetilde{B}^2_1+\widetilde C^2_1)(\widetilde{B}^2_2+\widetilde C^2_2)}\nonumber\\
   &=\left\{\int \Big|\frac tz\Big|+\frac{|tz|}{|z+c^{-1}_1tm(wm-1)|^2}dH_x(t)\int t+\frac{t}{|1+tc^{-1}_2w|^2}dH_y(t)\right.\nonumber\\
   &\quad\quad +\left.\int \frac{t^2c^{-1}_1|wm-1|^2}{|z+c^{-1}_1tm(wm-1)|^2}dH_x(t) \int \frac{t^2c^{-1}_2}{|1+tc^{-1}_2w|^2}dH_y(t)\right\}^{1/2}\nonumber\\
   &\quad \times \left\{\int \Big|\frac tz\Big|+\frac{|tz|}{|z+c^{-1}_1t\widetilde{m}(\widetilde{w}\widetilde{m}-1)|^2}dH_x(t)
      \int t+\frac{t}{|1+tc^{-1}_2\widetilde{w}|^2}dH_y(t)\right.\nonumber\\
   &\quad \quad  +\left.\int \frac{t^2c^{-1}_1|\widetilde{w}\widetilde{m}-1|^2}{|z+c^{-1}_1t\widetilde{m}(\widetilde{w}\widetilde{m}-1)|^2}dH_x(t)\int\frac{t^2c^{-1}_2}{|1+tc^{-1}_2\widetilde{w}|^2}dH_y(t)\right\}^{1/2}.
\end{align}
On the other hand, taking the imaginary part on both sides of the second equation in  \eqref{lsdt} and \eqref{lsdone}, we obtain
\begin{align}
  &\Im (\overline{m})=\int \frac{t^2c^{-1}_2\Im (w)}{|1+tc^{-1}_2w|^2}dH_y(t),\label{ima}\\
  &\Im (w)=\Im(wm\overline{z}-\overline{z})\int \frac{t}{|z|^2}+\frac{t}{|z+c^{-1}_1tm\big(wm-1\big)|^2}dH_x(t)\nonumber\\
  &\quad\quad \quad +\Im (\overline{m})\int \frac{t^2c^{-1}_1|wm-1|^2}{|z+c^{-1}_1tm(wm-1)|^2}dH_x(t).\label{imw}
\end{align}
Further, if it holds
\begin{align}
  \Im(wm\overline{z}-\overline{z})>|z|\Im(w)\int t+\frac{t}{|1+tc^{-1}_2w|^2}dH_y(t),\label{imwa}
\end{align}
then for $w \in \mathbb{C}^{+}$, combining the above three equations \eqref{ima}, \eqref{imw} and \eqref{imwa} will lead to
\begin{align}\label{ineqw}
  1&>\int \frac{t}{|z|}+\frac{t|z|}{|z+c^{-1}_1tm\big(wm-1\big)|^2}dH_x(t) \int t+\frac{t}{|1+tc^{-1}_2w|^2}dH_y(t)\nonumber\\
   &\quad \quad+\int \frac{t^2c^{-1}_1|wm-1|^2}{|z+c^{-1}_1tm(wm-1)|^2}dH_x(t) \int \frac{t^2c^{-1}_2}{|1+tc^{-1}_2w|^2}dH_y(t).
\end{align}
Such inequality also holds true if we replace $w$ and $m$ by $\tilde{w}$ and $\tilde{m}$, that is,
\begin{align}\label{ineqw2}
  1&>\int \frac{t}{|z|}+\frac{t|z|}{|z+c^{-1}_1t\tilde{m}\big(\tilde{w}\tilde{m}-1\big)|^2}dH_x(t) \int t+\frac{t}{|1+tc^{-1}_2\tilde{w}|^2}dH_y(t)\nonumber\\
   &\quad \quad+\int \frac{t^2c^{-1}_1|\tilde{w}\tilde{m}-1|^2}{|z+c^{-1}_1t\tilde{m}(\tilde{w}\tilde{m}-1)|^2}dH_x(t)\int \frac{t^2c^{-1}_2}{|1+tc^{-1}_2\tilde{w}|^2}dH_y(t).
\end{align}
Combining \eqref{ineqw} and \eqref{ineqw2} will lead to  a contradiction to \eqref{ine1}, which means that we could only have  one  solution $(s,w,m)$  satisfying  the system of equations \eqref{lsd-sys} on the set \eqref{set}.
\smallskip

So it is sufficient to prove the assertion \eqref{imwa} on some open set of $\mathbb C^+$. In fact,
using  the first and second equations in \eqref{lsdt}, we have
\begin{align*}
  &\Im(wm\overline{z}-\overline{z})=|z|^2 \Im(s),\\
  &\Im(zs)=\Im(wm)=\int t+\frac{t}{|1+tc^{-1}_2w|^2}dH_y(t)\Im(w).
\end{align*}
Then  assertion \eqref{imwa} is equivalent to
\begin{align}\label{ss}
  \Im (s)> \frac{1}{|z|} \Im(zs).
\end{align}
Actually, for any subsequence $\{n_k\}$ such that
$$
s_{n_k}(z)=\frac{1}{n_k}\E\tr(\B_{n_k}-zI_n)^{-1}
$$
converges,
the empirical distribution $F^{\B_{n_k}}$ has a limit $F$ (may depend on $\{n_k\}$), as $k\to\infty$, whose support is bounded upward by a constant, say $K$, which dose not depend on $\{n_k\}$. Moreover, the limit $s(z)$ of $s_{n_k}(z)$ is the Stieltjes transform of $F$, i.e.
\begin{align*}
  s(z)=\int \frac{1}{x-z}dF(x).
\end{align*}
This implies
\begin{align*}
  &\Im(s(z))=\int \frac{1}{|x-z|^2}dF(x) \Im(z),\\
  &\Im(zs(z))=\int \frac{x}{|x-z|^2}dF(x) \Im(z).
\end{align*}
Therefore, \eqref{ss} is true whenever $|z|> K$, which completes our proof.

\bigskip

\noindent{\em Non-Gaussian case:}
    %
since the two sets of samples $\{\xx_i\}$ and $\{\y_i\}$ are independent, we first fix the sequence of matrices $(\A_n)$ and show that, without the Gaussian assumption, the empirical spectral distribution $F^{\bS_{xy}}$ will still converge weakly to the same spectral distribution $F$ under Assumptions (a)-(c).
Next, the same trick can be applied to  $\{\xx_i\}$, which will not be detailed here. Our strategy to remove the Gaussian assumption is based on Lemma \ref{cth}, an extension of Lindeberg's argument for general smooth functions, see also Corollary 1.2 in \cite{ch08}. As a special case, letting $g$ be the identity function and $f$ be the Stieltjes transform, the theorem will ensure that the order of the difference in expectation between the two Stieltjes transforms under the Gaussian distribution and a non-Gaussian one is $O(n^{-1/2})$ whenever the two distributions match the first two moments and have finite fourth moment. Hence, such difference can be negligible as $n \to \infty$,  by which and the ``Step 1'' for Gaussian case the proof is done.  


Recall that
$$\B_n=\A^{1/2}_n\Big(\frac 1q \Y'\Y+\gamma_y \I\Big)\A^{1/2}_n=\A^{1/2}_n\Big(\frac 1q \W'\Sig_y \W+\gamma_y I\Big)\A^{1/2}_n,
$$
where the table $\W$ consists i.i.d.\ standard Gaussian random variables and we vectorize it as a $qn$-dimensional random vector, denoted as $\w=(w_{ij})$. Therefore, the Stieltjes transform $s_n(z)$ of $F^{\B_n}$ can be viewed as a function of the random vector $\w$, defined as
\begin{align*}
  U:= f(\w)=\frac 1n \tr (\B_n-z \I)^{-1},
\end{align*}
Similarly, we denote by
\begin{align*}
  V:= f(\tilde \w)
\end{align*}
the non-Gaussian counterpart of $U$, where
$\tilde \w=(\tilde w_{ij})$ have the same first two moments as $\{w_{ij}\}$ and finite fourth moment. Let $\bar \w=(\bar w_{ij})$ be a mixture of $\w$ and $\tilde \w$ by taking $\bar w_{ij}\in \{w_{ij}, \tilde w_{ij}\}$ for $i=p+1,\ldots,p+q$ and $j=1,\ldots,n$, whose matrix form is denoted by $\overbar W$.
Applying Lemma \ref{cth}, one gets
\begin{align}\label{u-v}
  |\E (U)-\E (V)|\leq K qn\lambda_3(f),
\end{align}
where
$$
\lambda_3(f)=\sup\left \{\bigg|\frac{\partial ^k f(\bar\w)}{\partial \bar w_{ij}^k}\bigg|^{3/k}: p+1\leq i\leq p+q, 1\leq j\leq n, 1\leq k\leq 3, \bar\w \in \mathbb R^{qn}\right\}.
$$
Hence, the remaining work is to find a bound for $\lambda_3(f)$, which can be achieved from bounding  the first three derivatives of $f$ with respect to $\bar w_{ij}$. To this end,
following the same truncation, centralization and rescaling steps as in  \cite{BSbook} (see Eq. (4.3.4)) and the ``no eigenvalues'' argument under finite fourth moment condition in \cite{BS98}, without loss of generality, we assume that the atoms $(\bar w_{ij})$ satisfy the following:
\begin{align*}
  &\E(\bar w_{ij})=0, ~\var(\bar w_{ij})=1,~|\bar w_{ij}|\leq \sqrt n,~ {\e_i'\overbar \W\overbar \W'\e_i\leq Kn,}
\end{align*}
for all $i$ and $j$,
where the vector {$\e_i$ is the $i$th canonical basis on $\mathbb R^q$}. For convenience, we still use notations $(w_{ij}, \w, \W)$ instead of $(\bar w_{ij}, \bar \w, \overbar \W)$ in what follows.

Let $\G=(\B_n-z I)^{-1}$, then
the first three derivatives of $f(\w)$ with respect to $ w_{ij}$ are the following:
\begin{align*}
  &\frac{\partial f(\w)}{\partial w_{ij}}=\frac 1n \tr \G'=-\frac 1n \tr \B'_n\G^2,\nonumber\\
  &\frac{\partial^2 f(\w)}{\partial w^2_{ij}}=-\frac 1n \tr (\B''_n\G^2+2\B'_n\G\G')=-\frac 1n \tr \B''_n\G^2+\frac 2n \tr \B'_n\G^2\B'_n\G\nonumber,\\
  &\frac{\partial^3 f(\w)}{\partial w^3_{ij}}=\frac 4n \tr \B''_n\G^2\B'_n-\frac 6n \tr \B'_n\G^2\B'_n\G\B'_n\G+\frac 2n \tr \B'_n\G^2\B''_n\G,
\end{align*}
where
\begin{align*}
  &\G'=-\G\B'_n\G,\\
  &\B'_n=\frac 1q \A^{1/2}_n(\e_j\e'_i\Sig_y\W+\W'\Sig_y\e_i\e'_j)\A^{1/2}_n,\nonumber\\
  &\B''_n=\frac 2q \A^{1/2}_n\e_j\e'_i\Sig_y\e_i\e'_j\A^{1/2}_n.
\end{align*}
and {the vector $\e_j$ is the $j$th canonical basis on $\mathbb R^n$}.

For the first derivative of $f$, since $\Sig_y$, $\A^{1/2}_n$ and $\G^2$ are all normal, we have
\begin{align}\label{ff}
  \sup\left|\frac{\partial f(\w)}{\partial w_{ij}}\right|&\leq \sup\left\{\frac {1}{nq} \left|\tr \A^{1/2}_n\e_j\e'_i\Sig_y\W\A^{1/2}_n\G^2\right|+\frac {1}{nq}\left| \tr \A^{1/2}_n\W'\Sig_y\e_i\e'_j\A^{1/2}_n\G^2\right|\right\}\nonumber\\
                                                         &\leq \sup\left\{\frac{K}{nq}\|\e_j\|\|\e'_i\W\|+\frac{K}{nq} \|\W'\e_i\|\|\e'_j\|\right\}\nonumber\\
                                                         &\leq Kn^{-3/2}.
\end{align}
For  the  second derivative, we have
\begin{align*}
  \left|\frac 1n \tr \B''_n\G^2\right|&=\frac {2}{nq} \left|\tr \A^{1/2}_n\e_j\e'_i\Sig_y\e_i\e'_j\A^{1/2}_n\G^2\right|
                                      \leq \frac{K}{nq}\|\e_j\|\cdot \|\e'_i\e_i\e'_j\|\leq Kn^{-2}~
\end{align*}
and
\begin{align*}
  &\quad~ \left|\frac 2n \tr \B'_n\G^2\B'_n\G\right|\nonumber\\
  &=\frac {2}{nq^2}\left| \tr \A^{1/2}_n(\e_j\e'_i\Sig_y\W+\W'\Sig_y\e_i\e'_j)\A^{1/2}_n\G^2 \A^{1/2}_n(\e_j\e'_i\Sig_y\W+\W'\Sig_y\e_i\e'_j)\A^{1/2}_n\G \right|\nonumber\\
  &\leq \frac{K}{nq^2}\left(\|\e_j\| \|\e'_i\W \e_j\e'_i\W\|+\|\e_j\| \|\e'_i\W\W'\e_i\e'_j\|+\|\W'\e_i\| \|\e'_j\e_j\e'_i\W\|+\|\W'\e_i\| \|\e'_j\W'\e_i\e'_j\|\right)\nonumber\\
  &\leq \frac{K}{nq^2}\left(n+\sqrt n \cdot |w_{ij}|\right)\nonumber\\
  &\leq Kn^{-2},
\end{align*}
which leads to the conclusion that
\begin{align}\label{sf}
  \sup\left|\frac{\partial^2 f(\w)}{\partial w^2_{ij}}\right|\leq Kn^{-2}.
\end{align}
Similarly, we could bound  the third derivative as follows,
\begin{align}\label{tf}
  \sup\left|\frac{\partial^3 f(\w)}{\partial w^3_{ij}}\right|&\leq \sup\left\{\frac{K}{nq^3}
                                                               \left(\|\e'_i\W\||w_{ij}|^2+2|w_{ij}||\e_i'\W\W'\e_i|+\|\e'_i\W\||\e_i'\W\W'\e_i|\right)\nonumber\right.\\
                                                             &\left.\quad +\frac{K}{nq^2}\left(\|\e'_i\W\|+|w_{ij}|\right)\right\}\nonumber\\
                                                             &\leq Kn^{-5/2}.
\end{align}
Finally,  combing \eqref{ff}, \eqref{sf} and \eqref{tf} gives
$$\lambda_3(f)=\sup \left\{\left|\frac{\partial f}{\partial w_{ij}}\right|^3, \left|\frac{\partial^2 f}{\partial w^2_{ij}}\right|^{\frac{3}{2}}, \left|\frac{\partial^3 f}{\partial w^3_{ij}}\right|\right\}=K n^{-5/2},$$ which together with \eqref{u-v} imply
\begin{align*}
  |\E (U)-\E (V)|\leq Kn^{-1/2}\to 0,\quad \text{as}\ n\to\infty.
\end{align*}
The proof is done.

\subsection{Proof of Theorem \ref{lsd-spike}}

Under our model setting \eqref{zx},  the three  data matrices $\X$, $\Y$ and $\Z$ are related as:
$$\Z=\Gamma \X\bS+\Y,$$ where $\Gamma=\sum_{k=1}^m\theta_k \bu_k \bv_k' $ and $\bS={\rm Diag}(\varepsilon_1, \ldots, \varepsilon_n)$.
So we have
\begin{align*}
  \frac 1q\Z'\Z&=\frac 1q \Y'\Y+\frac 1q\bS\X'\Gamma' \Gamma \X\bS+\frac 1q \bS\X'\Gamma'\Y+\frac 1q \Y'\Gamma \X\bS\nonumber\\
             &\triangleq\frac 1q \Y'\Y+\bH,
\end{align*}
where 
\begin{align}\label{zz}
\bH=\frac 1q\bS\X'\Gamma' \Gamma \X\bS+\frac 1q \bS\X'\Gamma'\Y+\frac 1q \Y'\Gamma \X\bS
\end{align} is a matrix of finite rank, at most $2m$.
Denote
$$
\widetilde \bS_{xz}=\A_n^{1/2}\left(\frac 1q \Z'\Z+\gamma_{z}\I_n\right)\A_n^{1/2}\quad\text{and}\quad
\widehat \bS_{xz}=\A_n^{1/2}\left(\frac 1q \Y'\Y+\gamma_{z}\I_n\right)\A_n^{1/2},\quad
$$
where
\begin{align}\label{kappa-yz}
 \gamma_{z}={\frac{1}{q}}\tr(\Sig_z)=\gamma_{y}+\frac{1}{q}\sum_{i=1}^m \theta^2_i \cdot \gamma_x=\gamma_{y}+o(1).
\end{align}
Applying Lemma \ref{rank-ineq} to $\B_n$, $\widetilde \bS_{xz}$ and $\widehat \bS_{xz}$, we have
\begin{align}\label{sss}
  ||F^{\widetilde \bS_{xz}}-F^{\widehat \bS_{xz}}||\to 0\quad\text{and}\quad L^3(F^{\B_n},F^{\widehat \bS_{xz}})\to 0,
\end{align}
almost surely, as $(n,p,q)$ tend to infinty. Combining \eqref{sss} and the fact that $\widetilde \bS_{xz}$ shares the same LSD as $\bS_{xz}$,
we conclude that $F^{\bS_{xz}}$ converges weakly to the LSD $F$  defined by \eqref{lsd}. The proof is thus complete.

\subsection{Proof of Theorem \ref{th-spike}}

We first note that, from the convergence in \eqref{k-r} and \eqref{kappa-yz}, asymptotically, the largest eigenvalues of $\bS_{xz}$ are the same as those of
$$
\bar \bS_{xz}:= \A^{1/2}_n\left(\frac 1q \Y'\Y+\bH+\gamma_{y}\I_n\right)\A^{1/2}_n,
$$
where $\bH$ is given in \eqref{zz}. So it's equivalent to prove the theorem for $\bar \bS_{xz}$.

\smallskip

\red{Next,  from \cite{BS98} and the inequality
$$
||\A^{1/2}_n\C_n\A^{1/2}_n||\leq ||\A_n||\cdot||\C_n||,
$$
we know that the spectral norm $||\A^{1/2}_n\C_n\A^{1/2}_n||$ is bounded in $n$, almost surely.
Define
$$
\lambda_+=\limsup_{n\to \infty} ||\A^{1/2}_n\C_n\A^{1/2}_n||,
$$
we consider the existence of spiked eigenvalues $(\lambda_{n, \ell})$ of $\bar \bS_{xz}$ in the interval $(\lambda_+,+\infty)$.
That is, for each $\ell \in \{1,\ldots, k\}$, $\lambda_{n, \ell}$ is an eigenvalue of $\bar \bS_{xz}$ but not an eigenvalue of  $\A^{1/2}_n\C_n\A^{1/2}_n$, i.e.
\begin{align}\label{ell-n}
  \left|\lambda \I_n-\bar \bS_{xz}\right|=0\quad\text{and}\quad \left|\lambda \I_n-\A^{1/2}_n\C_n\A^{1/2}_n\right|\neq 0,
\end{align}
for  $\lambda\in \{\lambda_{n,1},\ldots,\lambda_{n,k}\}$. 
}

In the following, we will show 
the  limits of $\lambda$ is defined in \eqref{outlier}.  
Under the assumptions in \eqref{ell-n}, we have
\begin{align}\label{ln-1}
  \left|\I_{n}-\left(\lambda \I_n-\A^{1/2}_n\C_n\A^{1/2}_n\right)^{-1} \A^{1/2}_n \bH\A^{1/2}_n\right|=0.
\end{align}
Recall the definition of $\bH$ in \eqref{zz},
then with  a little bit calculation, this matrix can be decomposed as
\begin{align}\label{h}
  \bH=\frac{1}{q}\left(\begin{array}{ccccc}
                            \ba_1 & \bb_1 & \cdots &\ba_m & \bb_m\end{array}\right)\left(\begin{array}{ccccc}
                                                              \theta_1 \lambda_{11} & 0  & \cdots & 0 & 0\\
                                                               0&\theta_1 \lambda_{12} & \cdots & 0 & 0\\
                                                               \vdots &  \vdots& \ddots &  \vdots& \vdots\\
                                                               0& 0 & \cdots & \theta_m\lambda_{m1}& 0\\
                                                               0& 0 & \cdots &0& \theta_m \lambda_{m2}
                                                             \end{array}\right)\left(\begin{array}{c}
                                                                                       \ba'_1\\
                                                                                       \bb'_1\\
                                                                                       \vdots\\
                                                                                       \ba'_m\\
                                                                                       \bb'_m\\
                                                                                     \end{array}\right),
\end{align}
where
\begin{align*}
  \ba_i&=u_{i1}\bS\X'\bv_i+w_{i1}\Y'\bu_i,\nonumber\\
  \bb_i&=u_{i2}\bS\X'\bv_i+w_{i2}\Y'\bu_i,\nonumber\\
  \lambda_{i1}&=\|\bS\X'\bv_i\|\|\Y'\bu_i\|\left\{\frac{\sqrt{4\|\Y'\bu_i\|^2+\theta^2_i\|\bS\X'\bv_i\|^2}+\theta_i \|\bS\X'\bv_i\| }{\sqrt{4\|\Y'\bu_i\|^2+\theta^2_i\|\bS\X'\bv_i\|^2}-\theta_i \|\bS\X'\bv_i\|}\right\}^{1/2},\nonumber\\
  \lambda_{i2}&=-\|\bS\X'\bv_i\|\|\Y'\bu_i\|\left\{\frac{\sqrt{4\|\Y'\bu_i\|^2+\theta^2_i\|\bS\X'\bv_i\|^2}-\theta_i \|\bS\X'\bv_i\| }{\sqrt{4\|\Y'\bu_i\|^2+\theta^2_i\|\bS\X'\bv_i\|^2}+\theta_i \|\bS\X'\bv_i\|}\right\}^{1/2},
\end{align*}
with
\begin{align*}
  &u_{i1}=\frac{1}{\|\bS\X'\bv_i\|}\left\{\frac 12+\frac{\theta_i \|\bS\X'\bv_i\|}{2\sqrt{4\|\Y'\bu_i\|^2+\theta^2_i\|\bS\X'\bv_i\|^2}}\right\}^{1/2},\nonumber\\
  &u_{i2}=\frac{1}{\|\bS\X'\bv_i\|}\left\{\frac 12-\frac{\theta_i \|\bS\X'\bv_i\|}{2\sqrt{4\|\Y'\bu_i\|^2+\theta^2_1\|\bS\X'\bv_i\|^2}}\right\}^{1/2},\nonumber\\
  &w_{i1}=\frac{1}{\|\Y'\bu_i\|}\left\{\frac 12-\frac{\theta_i \|\bS\X'\bv_i\|}{2\sqrt{4\|\Y'\bu_i\|^2+\theta^2_i\|\bS\X'\bv_i\|^2}}\right\}^{1/2},\nonumber\\
  &w_{i2}=-\frac{1}{\|\Y'\bu_i\|}\left\{\frac 12+\frac{\theta_i \|\bS\X'\bv_i\|}{2\sqrt{4\|\Y'\bu_i\|^2+\theta^2_i\|\bS\X'\bv_i\|^2}}\right\}^{1/2}.
\end{align*}
In addition, it's straightforward to verify the following relations,
\begin{align}\label{luw}
  \left\{\begin{array}{l}
         \displaystyle  \lambda_{i1}u^2_{i1}+\lambda_{i2}u^2_{i2}=\theta_i, \\
         \displaystyle   \lambda_{i1}w^2_{i1}+\lambda_{i2}w^2_{i2}=0,\\
         \displaystyle   \lambda_{i1}u_{i1}w_{i1}+\lambda_{i2}u_{i2}w_{i2}=1.
         \end{array}\right.
\end{align}
Denote $\D_n= \A^{1/2}_n\left(\lambda \I_n-\A^{1/2}_n\C_n\A^{1/2}_n\right)^{-1}\A^{1/2}_n$ and
\begin{align*}
\M_n=  \frac{1}{q}\left(\begin{array}{c}
                                                                                       \ba'_1\\
                                                                                       \bb'_1\\
                                                                                       \vdots\\
                                                                                       \ba'_m\\
                                                                                       \bb'_m\\
                                                                                     \end{array}\right) \D_n\left(\begin{array}{ccccc}
                            \ba_1 & \bb_1 & \cdots &\ba_m & \bb_m\end{array}\right)\left(\begin{array}{ccccc}
                                                              \theta_1 \lambda_{11} & 0  & \cdots & 0 & 0\\
                                                               0&\theta_1 \lambda_{12} & \cdots & 0 & 0\\
                                                               \vdots &  \vdots& \ddots &  \vdots& \vdots\\
                                                               0& 0 & \cdots & \theta_m\lambda_{m1}& 0\\
                                                               0& 0 & \cdots &0& \theta_m \lambda_{m2}
                                                             \end{array}\right).
\end{align*}
Then  \eqref{ln-1}  and  \eqref{h} imply
\begin{align*}
f_n(\lambda) := |\I_{2m}-\M_n|=0.
\end{align*}

We next find the limit of $f_n(\lambda)$.
Let
\begin{align*}
  \alpha_{n}=\frac 1n\tr \bS\D_n\bS(\A_n-\gamma_x \I_n)\quad \text{and}\quad
  \beta_{n}=\frac 1n\tr \D_n(\C_n-\gamma_y \I_n),
\end{align*}
one may get for any $i \in \{1,\ldots, m\}$,
\begin{align*}
 & \frac { \ba'_i\D_n\ba_i}{q}=\frac{u_{i1}^2}{c_{n2}}\alpha_{n}+\frac{w_{i1}^2}{c_{n2}}\beta_{n}+o_{a.s.}(1)
  ,\\
&\frac { \ba'_i\D_n\bb_i}{q}=\frac{u_{i1}u_{i2}}{c_{n2}}\alpha_{n}+\frac{w_{i1}w_{i2}}{c_{n2}}\beta_{n}+o_{a.s.}(1)
  ,\\
 & \frac{
  \bb'_i\D_n\bb_i}{q}=\frac{u_{i2}^2}{c_{n2}}\alpha_{n}+\frac{w_{i2}^2}{c_{n2}}\beta_{n}+o_{a.s.}(1)~,
\end{align*}
and for any  $i \neq j \in \{1,\ldots, m\}$,
\begin{align*}
 \frac { \ba'_i\D_n\ba_j}{q}=o_{a.s.}(1)~,\quad  \frac { \ba'_i\D_n\bb_j}{q}=o_{a.s.}(1)~.
 \end{align*}
>From the above approximations and the identities in \eqref{luw}, we have

\begin{align*}
f_n(\lambda)=\prod_{k=1}^m\big|\I_2-\M_{nk}\big|+o_{a.s}(1)
\end{align*}
where
\begin{align}\label{mnk}
\M_{nk}=\frac{\theta_k}{c_{n2}}
\left(\begin{array}{cc}
\alpha_n  &  0\\
0 & \beta_n
\end{array}
\right)
\left(\begin{array}{cc}
\theta_k & 1\\
1 & 0
\end{array}
\right).
\end{align}
Let
$\be=(\varepsilon_1, \ldots,\varepsilon_n)'$, then
\begin{align}\label{eqq1}
\alpha_n=\frac 1n\tr \bS\D_n\bS\A_n-\frac{\gamma_x}{n}\tr\bS\D_n\bS=\frac 1n \be'\left(\D_n\circ \A_n\right) \be-\frac{\gamma_x}{n}\be' \text{Diag}(\D_n)\be,
\end{align}
where ``$\circ$'' denotes the Hadamard product of two matrices.
According to Theorem 1 of \cite{v68}, we have
\begin{align}\label{eqq2}
&\frac 1n \be'\left(\D_n\circ \A_n\right) \be-\frac 1n\E  \left[\be'\left(\D_n\circ \A_n\right) \be \right]\xrightarrow{a.s.} 0,\\
&\frac{1}{n}\be' \text{Diag}(\D_n)\be-\frac 1n \E \tr \D_n\xrightarrow{a.s.} 0.
\end{align}
Further, 
\begin{align}\label{eqq3}
&~\quad \frac 1n\E  \left[\be'\left(\D_n\circ \A_n\right) \be\right]=\frac 1n \E \tr \big[\D_n \text{Diag}(\A_n)\big]\nonumber\\
&=\frac 1n \E \tr \big[\D_n \left(\text{Diag}(\A_n)-2\gamma_x\I_n\right)\big]+\frac {2\gamma_x}{n} \E \tr \D_n\nonumber\\
&=\frac {2\gamma_x}{n} \E \tr \D_n+o(1),
\end{align}
where the last equality is due to the following convergence,
\begin{align*}
\left|\frac 1n \tr \big[\D_n\cdot \left(\text{Diag}(\A_n)-2\gamma_x\I_n\right)\big]\right|\leq \frac 1n \|\D_n\|\cdot\tr\big|\A_n-2\gamma_x\I_n\big|\xrightarrow{a.s.} 0.
\end{align*}
Collecting results in \eqref{eqq1}-\eqref{eqq3}, we get
\begin{align}\label{alpha-n}
\alpha_n=-\gamma_xw_n(\lambda)+o_{a.s.}(1)\xrightarrow{a.s.}\alpha\triangleq -w(\lambda)\int tdH_x(t),
\end{align}
where $w_n(z)$ is defined in \eqref{wmn}, whose domain can be expanded to $(\lambda_+, +\infty)$ for all large $n$.
For $\beta_n$, we have
\begin{align}
  \beta_{n}&=\frac 1n\tr (\D_n\C_n)-\frac{\gamma_y}{n} \tr \D_n\no
  &=-1+\frac{\lambda}{n}\tr \left(\lambda \I_n-\A^{1/2}_n\C_n\A^{1/2}_n\right)^{-1}-\frac{\gamma_y}{n} \tr \D_n\no
  &=-\frac{1}{n}\sum_{k=1}^q\frac{\tau_kw_{n}(\lambda)}{c_{n2}+\tau_kw_{n}(\lambda)}+o_{a.s.}(1)\nonumber\\
           &\xrightarrow{a.s.}\beta\triangleq -c_2\int\frac{tw(\lambda)dH_y(t)}{c_2+tw(\lambda)},\label{beta-n}
\end{align}
where the third equality is from \eqref{lsd:c2} with $(\tau_k)$ being the eigenvalues of $\Sig_y$. Collecting results in \eqref{mnk},\eqref{alpha-n} and \eqref{beta-n}, we get
$$
f_n(\lambda)\xrightarrow{a.s.}f(\lambda)\triangleq \prod_{k=1}^m\left(1-\theta_k^2g(\lambda)\right),
$$
where the function $g$ is given in \eqref{function-g}.
  With the definition of the critical value $\theta_0$ in
  \eqref{theta0}, we find that
 for any $k\in \{1,\ldots,m\}$ and $\theta_k>\theta_0$, there are $k$ zeros $\lambda_1>\cdots>\lambda_k$ of $f(\lambda)$ on $(\lambda_+,\infty)$.
By continuity arguments, see Lemma 6.1 in \cite{BN11}, we verify the existence of the spikes $\lambda_{n,1},\ldots,\lambda_{n,k}$ whose limits are $\lambda_1,\ldots,\lambda_{k}$, respectively.
The proof is then complete.

\end{document}